\documentstyle{amsppt}
\input epsf
\expandafter\let\csname logo\string @\endcsname=\empty
\define\a{\alpha}
\predefine\barunder{\b}
\redefine\b{\beta}
\define\bpage{\bigpagebreak}

\define\cd{\cdot}
\predefine\dotunder{\d}
\redefine\d{\delta}
\define\e{\epsilon}

\define\f{\frac}
\define\g{\gamma}
\define\G{\Gamma}
\define\k{kappa}
\define\lb\{{\left\{}

\define\lm{\limits}
\define\lra{\longrightarrow}

\define\Om{\Omega}
\define\om{\omega}

\define\oper{\operatorname}

\define\ov{\overline}

\define\p{\partial}
\define\rb\}{\right\}}

\define\sub{\subheading}

\define\th{\theta}

\define\vp{\varphi}
\define\wh{\widehat}

\define\x{\times}
\define\z{\zeta}
\define\({\left(}
\define\){\right)}
\define\[{\left[}
\define\]{\right]}
\define\<{\left<}
\define\>{\right>}
\def\slantline#1#2#3#4#5{\hbox to 0pt{
}}

%SCRIPT LETTERS
%%%%%%%%%%

\def\SF{{\Cal F}}

\def\SK{{\Cal K}}

\def\SR{{\Cal R}}
\def\SS{{\Cal S}}

\def\SY{{\Cal Y}}
\def\SZ{{\Cal Z}}

%%%%%%%
%%%BOLDFACE

\def\BC{{\bold C}}

\def\BE{{\bold E}}

\def\BH{{\bold H}}

\def\BR{{\bold R}}

%%%%%BLACKBOARD BOLD LETTERS

\def\bbc{{\Bbb C}}

\def\bbr{{\Bbb R}}

\magnification=1200
\pagewidth{5.4in}
\pageheight{7.5in}
\NoRunningHeads

\define\ext{\oper{Ext}}
\define\grad{\oper{grad}}
\define\id{\oper{id}}
\define\imm{\oper{Im}}
\define\nes{\text{NE}}
\define\qd{\oper{QD}}
\define\ree{\oper{Re}}
\define\sw{\text{SW}}
\define\tecs{Teichm\"uller's\ }
\define\Ree{\oper{Re}}
\define\im{\oper{Im}}
\define\supp{\oper{supp}}
\define\tec{Teich\-m\"uller\ }
\define\wei{Weier\-stra{\ss}\ }

\define\fin{finite\ }
\def\nez#1{#1_{\nes}}
\def\phine{\phi^{\nes}_k}
\def\phisw{\phi^{\sw}_k}
\def\psine#1{\psi^{\nes}_#1}
\def\psisw#1{\psi^{\sw}_#1}
\def\sww#1{#1_{\sw}}
\def\tsym#1{#1^{\text{symm}}_{0,2p+2}}

\def\sym{{\text{symm}}}
\def\tp{T_{0,2p+2}}

\def\Ts{T}

\topmatter
\title Minimal Surfaces of Least Total Curvature and Moduli
Spaces of Plane Polygonal Arcs
\endtitle
\affil  Matthias Weber\\
Mathematisches Institut\\
Universit\"at Bonn\\
Beringstra{\ss}e 6\\
Bonn, Germany\\
\\
Michael Wolf*\\
Dept. of Mathematics \\
Rice University \\ Houston, TX 77251 USA
\endaffil
\thanks *Partially supported by NSF grant number DMS 9300001 and the
Max Planck Institut; Alfred P\. Sloan Research Fellow
\endthanks
\endtopmatter

\document

\centerline{Revised April 13, 1998}

\sub{Introduction} The first major goal of this paper is to prove the
existence of complete minimal surfaces of each genus $p>1$ which minimize the
total curvature
(equivalently, the degree of the Gau{\ss} map)
for their genus. The genus zero version of these
surfaces is known as Enneper's surface (see \cite{Oss2}) and the genus
one version is due to Chen-Gackstatter (\cite{CG}). Recently,
experimental evidence for the existence of these surfaces for genus
$p\le35$ was found by Thayer (\cite{Tha}); his surfaces, like those
in this paper, are hyperelliptic surfaces with a single end, which is
asymptotic to the end of Enneper's surface.

Our methods for constructing these surfaces are somewhat novel, and
as their development is the second major goal of this paper, we
sketch them quickly here. As in the construction of other recent
examples of complete immersed (or even embedded) minimal surfaces in
$\BE^3$, our strategy centers around the \wei representation
for minimal surfaces in space, which gives a parametrization of the
minimal surface in terms of meromorphic data on the Riemann surface
which determine three meromorphic one-forms on the
underlying Riemann surface.

The art in finding a minimal surface via this representation lies in
finding a Riemann surface and meromorphic data on that surface so
that the representation is well-defined, i.e., the local \wei
representation can be continued around closed curves without
changing its definition. This latter condition amounts to a
condition on the imaginary parts of some periods of forms associated
to the original \wei data.

In many of the recent constructions of complete minimal
surfaces, the geometry of the desired surface is used to set up a space
of possible \wei data and Riemann surfaces, and then to consider the
period
problem as a purely analytical one. This approach is very effective
as long as the dimension of the space of candidates remains small.
This happens for instance if
enough symmetry of the resulting surface is assumed so that
the moduli space of candidate (possibly singular) surfaces has very
small dimension -- in fact, there is sometimes only a single surface to
consider.  Moreover, the candidate quotient surfaces (for instance,
a thrice punctured sphere) often have a relatively well-understood
function theory which serves to simplify the space of possibilities,
even if the (quite difficult) period problem for the
\wei data still remains.

In our situation however, the dimension of the space
of candidates grows with the genus.
Our approach is to first view
the periods and the conditions on them
as defining
a geometric object (and inducing a construction of a
pair of Riemann surfaces), and to then
prove analytically that the Riemann surfaces are identical, employing
methods from \tec theory.

Generally speaking, our approach is to
construct two different Riemann surfaces, each with a meromorphic
one-form, so that the period problem would be solved if only the
surfaces would coincide. To arrange for a situation where we can
simultaneously define a Riemann surface, and a meromorphic one-form
on that surface with prescribed periods, we exploit the perspective
of a meromorphic one-form as defining a singular flat structure on
the Riemann surface, which we can develop onto $\BE^2$.

In particular, we first assume sufficient symmetry of the Riemann
surface so that the quotient orbifold flat structure has a
fundamental domain in $\BE^2$ which is bounded by a properly embedded
arc composed of $2p+2$ horizontal and vertical line segments with
the additional properties that the segments alternate from horizontal
segments to vertical segments, with the direction of travel also
alternating between left and right turns. We call such an arc a
'zigzag'; further we restrict our attention to `symmetric zigzags',
those zigzags which are symmetric about the line $\{y=x\}$ (see
Figure~1).

A crucial observation is that we can turn this construction around.
Observe that a zigzag $Z$ bounds two domains, one, $\nez\Om(Z)$,
on the northeast side, and one, $\sww\Om(Z)$, on the southwest
side. When we double each of these domains and then take a double
cover of the resulting surface, branched over each of the images of
the vertices of $Z$, we have two hyperelliptic Riemann surfaces,
$\nez\SR(Z)$ and $\sww\SR(Z)$, respectively. Moreover, the form
$dz$ when restricted to $\nez\Om(Z)$ and $\sww\Om(Z)$, lifts to
meromorphic one-forms $\nez\om(Z)$ on $\nez\Om(Z)$ and
$\sww\om(Z)$ on $\sww\Om(Z)$ both of whose sets of periods are
integral linear combinations of the periods of $dz$ along the
horizontal and vertical arcs of $Z$.

Then, suppose for a moment that we can find a zigzag $Z$ so that
$\nez\Om(Z)$ is conformally equivalent to $\sww\Om(Z)$ with the
conformal equivalence taking vertices to vertices (where $\infty$ is
considered a vertex). (We call such
a zigzag {\it reflexive}.)  Then $\nez\SR(Z)$ would be conformally
equivalent to $\sww\SR(Z)$ in a way that $e^{i\pi/4}\nez\om(Z)$
and $e^{-i\pi/4}\sww\om(Z)$ have conjugate periods. As these forms
will represent $gdh$ and $g^{-1}dh$ in the classical \wei
representation $X=\ree\int(\f12\(g-\f1g\)$, $\f i2\(g+\f1g\)$, 1) $dh$,
it will turn out that this conformal equivalence is just what we need
for the \wei representation based on $\nez\om(Z)$ and $\sww\om(Z)$
to be well-defined.

This construction is described precisely in \S3.

We are left to find such a zigzag. Our approach is non-constructive
in that we consider the space $\SZ_p$ of all possible symmetric
zigzags with $2p+2$ vertices and then seek, within that space $\SZ$,
a symmetric zigzag $Z_0$ for which there is a conformal equivalence
between $\nez\Om(Z_0)$ and $\sww\Om(Z_0)$ which preserves vertices.
The bulk of the paper, then, is an analysis of this moduli
space $\SZ_p$ and some functions on it, with the goal of finding a
certain fixed point within it.

Our methods, at least in outline, for finding such a symmetric zigzag
are quite standard in contemporary \tec theory. We first find that
the space $\SZ_p$ is topologically a cell, and then we seek an
appropriate height function on it. This appropriate height function should be
proper, so that it has an interior critical point, and it should have
the feature that at its critical point $Z_0\in\SZ_p$, we have the
desired vertex-preserving conformal equivalence between
$\nez\Om(Z_0)$ and $\sww\Om(Z_0)$.

One could imagine that a natural height function might be the \tec
distance between $\nez\SR(Z)$ and $\sww\SR(Z)$, but it is easy to
see that there is a family $Z_t$ of zigzags, some of whose vertices
are coalescing, so that the \tec distance between $\nez\SR(Z_t)$ and
$\sww\SR(Z_t)$ tends to a finite number. We thus employ a different
height function $D(\cd)$ that, in effect, blows up small scale
differences between $\nez\SR(Z_t)$ and $\sww\SR(Z_t)$ for a
family of zigzags $\{Z_t\}$ that leave all compacta of $\SZ_p$.

We discuss the space $\SZ_p$ of symmetric zigzags, and study
degeneration in that space in \S4. We show that
the map between the marked extremal length spectra for $\nez\SR(Z)$ and
$\sww\SR(Z)$ is not real analytic at infinity in $\SZ_p$, and thus there must
be small scale differences between those extremal length spectra.
We do this by first observing that both extremal lengths and
Schwarz-Christoffel integrals can be computed using generalized
hypergeometric functions; we then show that the well-known
monodromy properties of these functions lead to a crucial sign
difference in the asymptotic
expansions of the Schwarz-Christoffel maps at regular singular points.
Finally, these sign differences are
exploited to yield the desired non-analyticity.

Our height function $D(Z)$ while not the \tec distance between
$\nez\SR(Z)$ and $\sww\SR(Z)$, is still based on differences
between extremal lengths on those surfaces, and in effect, we follow
the gradient flow $dD$ on $\SZ_p$
from a convenient initial point in $\SZ_p$
to a solution of our problem.
There are two aspects to this approach.  First, it
is especially convenient that we know a formula (\cite{Gar}) for
$d[\ext_{[\g]}(\SR)]$ where $\ext_{[\g]}(\SR)$ denotes the extremal
length of the curve family $[\g]$ on a given Riemann surface $\SR$. This
gradient of extremal length is given in terms of a holomorphic
quadratic differential $2\Phi_{[\g]}(\SR)=d\ext_{[\g]}(\SR)$ and can
be understood in terms of the horizontal measured foliation of that
differential providing a `direction field' on $\SR$ along which to
infinitesimally deform $\SR$ as to infinitesimally increase
$\ext_{[\g]}(\SR)$. We then show that grad~$D(\cd)\bigm|_Z$ can be
understood in terms of a pair of holomorphic quadratic differentials
on $\nez\SR(Z)$ and $\sww\SR(Z)$, respectively, whose (projective)
measured foliations descend to a well-defined projective class of
measured foliations on $\BC=\nez\Om(Z)\cup Z\cup\sww\Om(Z)$. This
foliation class then indicates a direction in which to
infinitesimally deform $Z$ so as to infinitesimally decrease
$D(\cd)$, as long as $Z$ is not critical for $D(\cd)$. Thus, a
minimum for $D(\cd)$ is a symmetric zigzag $Z_0$ for which
$\nez\Om(Z_0)$ is conformally equivalent to $\sww\Om(Z)$ in a
vertex preserving way.  Second, it is technically convenient
to flow along a path in $\SZ_p$ in which the form of the height
function simplifies.  In fact, we will flow from a genus $p-1$ solution
in $\SZ_{p-1} \subset \p\ov\SZ$ along a path $\SY$ in $\SZ_p$ to a
genus $p$ solution
in $\SZ_p$.  Here the technicalities are that some \tec theory
and the symmetry we have imposed on the zigzags allow us to invoke
the implicit function theorem at a genus $p-1$ solution in $\p\ov\SZ$
to find such a good path $\SY \subset \SZ_p$.
Thus, formally,  we find a solution in each genus
inductively,  by showing that given a reflexive zigzag of genus
$p-1$, we can 'add a handle' to obtain a solution of genus $p$.

We study the gradient flow of the height function $D(\cd)$ in \S5.

Combining the results in sections 4 and 5, we conclude

\proclaim{Main Theorem B} There exists a reflexive symmetric zigzag
of genus ~$p$ for $p\ge0$ which is isolated in $\SZ_p$.
\endproclaim

When we interpret this result about zigzags as a result on \wei
data for minimally immersed Riemann surfaces in $\BE^3$, we derive as
a corollary

\proclaim{Main Theorem A} For each $p\ge 0$, there exists a
minimally immersed Riemann surfaces in $\BE^3$ with one Enneper-type
end and total curvature $-4\pi(p+1)$. This surface has at most eight
self-isometries.
\endproclaim

In \S6, we adapt our methods slightly to prove the existence of minimally
immersed surfaces of genus $p(k-1)$ with one Enneper-type end of
winding order $2k-1$: these surfaces extend and generalize examples
of Karcher ([Kar]) and Thayer ([Tha]), as well as those constructed
in Theorem A.

While we were preparing this manuscript several years ago
we received a copy of a
preprint by K. Sato [S] which also asserts Theorem A.  Our approach
is different than that of Sato, and possibly more general, as it is
possible to assign zigzag configurations to a number of families
of putative minimal surfaces. We discuss further applications of this
technique in a forthcoming paper [WW].

The authors wish to thank Hermann Karcher for many hours of pleasant
advice.

\newpage
\sub{\S2. Background and Notation}

\sub{2.1 Minimal Surfaces and the \wei Representation}Here
we recall some well-known facts from the theory of minimal
surfaces and put our result into context.

Locally, a minimal surface can always be described by {\it \wei
data}, i.e. there are always a simply connected domain $U$, a
holomorphic function $g$ and a holomorphic $1$-form $dh$ in $U$ such
that the minimal surface is locally given by

$$
z \mapsto \ree \int_\cd^z
\matrix
\f12(g-\frac1 g) dh \\
\frac i 2(g+\frac 1 g) dh\\
dh
\endmatrix
$$
For instance, $g(z)=z$ and $dh = \frac{dz}z$ will lead to the
catenoid, while $g(z)=z$ and $dh = zdz$ yields the Enneper surface.

It is by no means clear how global properties of a minimal surface
are related to this local representation.

However, two global properties together have very strong
consequences on the \wei data. One is the metrical completeness, and
the other the total (absolute) Gau{\ss}ian curvature of the surface
$R$, defined by
$$
\SK := \int_R \vert K \vert dA = 4\pi \cd \hbox{degree of the
Gau{\ss} map}
$$

We will call a complete minimal surface of finite absolute
Gau{\ss}ian curvature a {\it \fin} minimal surface.

Then by a famous theorem of R\. Osserman, every \fin minimal surface
(see \cite{Oss1, Oss2, Laws}) can be represented by \wei data which
are defined on a compact Riemann surface $R$, punctured at a finite
number of points. Furthermore, the \wei data extend to meromorphic
data on the compact surface. Thus, the construction of such surfaces
is reduced to finding meromorphic \wei data on a compact Riemann
surface such that the above representation is well defined, i.e.
such that all three $1$-forms showing up there have purely imaginary
periods. This is still not a simple problem.

\smallpagebreak

From now on, we will restrict our attention to \fin minimal surfaces.

Looked at from far away, the most visible parts of a \fin minimal
surface will be the ends. These can be seen from the \wei data by
looking at the singularities $P_j$  of the Riemannian metric which
is given by the formula
$$
ds = \(\vert g\vert + \frac 1 {\vert g\vert}\) \vert dh\vert
\tag2.0$$
An end occurs at a puncture $P_j$ if and only if $ds$ becomes
infinite in the compactified surface at $P_j$.

To each end is associated its winding or spinning number $d_j$, which
can be defined geometrically by looking at the intersection curve of
the end with a very large sphere which will be close to a great
circle and by taking its winding number, see \cite{Gack, J-M}. This
winding number is always odd: it is $1$ for the catenoid end $3$ for
the Enneper end. There is one other end of winding number $1$,
namely the planar end which of course occurs as the end of the
plane, but also as one end of the Costa surface, see \cite{Cos1}.
For \fin minimal surfaces, there is a Gau{\ss}-Bonnet formula
relating the total curvature to genus and winding numbers:
$$
\int_R K dA = 2\pi\(2(1-p)-r-\sum_{j=1}^r d_j\)
$$
where $p$ is the genus of the surface, $r$ the number of ends and
$d_j$ the winding number of an end. For a proof, see again
\cite{Gack, J-M}.

From this formula one can conclude that for a non planar surface
$\int_R \vert K\vert dA\le 4\pi(p+1)$ which raises the question of
finding for each genus a (non-planar)
minimal surface for which equality holds.
This is the main goal of the paper.

The following is known:
\roster
\item"$p=0$:"
The Enneper surface and the catenoid are the only non-planar minimal surfaces
with  $\SK=4\pi$, see e.g\. \cite{Oss2}.
\item"$p=1$:"
The only surface with $\SK=8\pi$ is the Chen-Gackstatter surface
(which is defined on the square torus), see \cite{CG, Lop, Blo}.
\item"$p=2$:"  An example with $\SK=12\pi$ was also constructed in
\cite{CG}. Uniqueness is not known here.
\item"$p=3$:"  An example with $\SK=16\pi$ was constructed by do
Esp\'{\i}rito-Santo (\cite{Esp}).
\item"$p\le35$:" E\. Thayer has solved the period problem
numerically and produced pictures of surfaces with minimal $\SK$.
\endroster

Note that all these surfaces are necessarily not embedded: For given
genus, a finite minimal surface of minimal $\SK$ could, by the winding
number formula, have only either one end of Enneper type (winding
number $3$) which is not embedded or two ends of winding number $1$.
But by a theorem of R\. Schoen (\cite{Sch}), an embedded finite
minimal surface with only two ends has to be the catenoid. Hence if
one looks for embedded minimal surfaces, one has to allow more $\SK$.
For the state of the art here, see \cite{Ho-Ka}.

If one allows even more total curvature and permits non-embeddedness,
some general methods are available, as explained in \cite{Kar}.

\sub{2.2. Zigzags} A {\it zigzag} $Z$ of genus $p$ is an open and properly
embedded arc in $C$ composed of alternating
horizontal and vertical subarcs with
angles of $\pi/2$, $3\pi/2$, $\pi/2$, $3\pi/2,\dots,\pi/2$ between
consecutive sides, and having $2p+1$ vertices ($2p+2$ sides,
including an initial infinite vertical side and a terminal infinite
horizontal side.) A {\it symmetric\/} zigzag of genus $p$ is a zigzag
of genus $p$ which is symmetric about the line $\{y=x\}$. The space
$\SZ_p$ of genus $p$ zigzags consists of all symmetric zigzags of
genus $p$
up to similarity; it is equipped with the topology induced by the embedding
of $\SZ_p\lra\BR^{2p}$ which associates to a zigzag $Z$ the $2p$-tuple
of its lengths of sides, in the natural order.

A symmetric zigzag $Z$ divides the plane $\BC$ into two regions, one
which we will denote by $\nez\Om(Z)$ which contains large positive
values of $\{y=x\}$, and the other which we will denote by
$\sww\Om(Z)$.  (See Figure 1.)

\epsfxsize 4in
\centerline{\epsfbox{figure1.ps}}
\centerline{Figure 1}

\sub{Definition 2.2.1} A symmetric zigzag $Z$ is
called {\it reflexive\/} if there is a
conformal map $\phi:\nez\Om(Z)\to\sww\Om(Z)$ which takes vertices to
vertices.

\example{Examples 2.2.2}
There is only one zigzag of genus $0$, consisting of the positive
imaginary and positive real half-axes. It is automatically symmetric
and reflexive.

Every symmetric zigzag of genus $1$ is also automatically reflexive.
\endexample

%\vskip 0.2in

\sub{2.3. \tec Theory} For $M$ a smooth surface, let Teich~$(M)$
denote the \tec space of all conformal structures on $M$ under the
equivalence relation given by pullback by diffeomorphisms isotopic
to the identity map id:~$M\lra M$. Then it is well-known that
Teich~$(M)$ is a smooth finite dimensional manifold if $M$ is a
closed surface.

There are two spaces of tensors on a Riemann surface $\SR$ that are
important for the \tec theory. The first is the space QD$(\SR)$ of
holomorphic quadratic differentials, i.e., tensors which have the
local form $\Phi=\vp(z)dz^2$ where $\vp(z)$ is holomorphic. The
second is the space of Beltrami differentials Belt$(\SR)$, i.e.,
tensors which have the local form $\mu=\mu(z)d\bar z/dz$.

The cotangent space $T^*_{[\SR]}$(Teich~$(M)$) is canonically
isomorphic to QD$(\SR)$, and the tangent space is given by
equivalence classes of (infinitesimal) Beltrami differentials, where
$\mu_1$ is equivalent to $\mu_2$ if
$$
\int_\SR\Phi(\mu_1-\mu_2) = 0\qquad\text{for every }\
\Phi\in\text{QD}(\SR).
$$

If $f:\BC \to \BC$ is a diffeomorphism, then the Beltrami
differential associated to the pullback conformal structure
is $\nu = \frac{f_{\bar z}}{f_z} \frac{d\bar z}{dz}$. If $f_{\e}$ is a
family of such diffeomorphisms with $f_0= identity$, then the
infinitesimal Beltrami differential is given by
$\frac d{d\e}\bigm|_{\e=0} \nu_{f_\e} 
=\(\frac d{d\e}\bigm|_{\e=0} f_\e\)_{\bar z}$.  
We will
carry out an example of this computation in \S5.2.

A holomorphic quadratic differential comes with a picture that is
a useful aid to one's intuition about them. The picture is that of
a pair of transverse measured foliations, whose properties we sketch
briefly (see \cite{FLP}, \cite{Ke}, and \cite{Gar} for more details). 
We next define a measured foliation on a (possibly) punctured Riemann
surface; to set notation, in what follows,
the Riemann surface $\SR$ is possibly punctured, i.e. there is
a closed Riemann surface $\overline{\SR}$, and a set of 
points $\{q_1, \dots, q_m\}$, so that 
$\SR = \overline{\SR} - \{q_1, \dots, q_m\}$.

A measured foliation $(\SF, \mu)$ on a Riemann surface $\SR$ with
singularities $\{p_1, \dots, p_l\}$ (where some of
the singularities might also be elements of the
puncture set $\{q_1, \dots, q_m\}$) consists of a foliation
$\SF$ of $\SR -\{p_1, \dots, p_l\}$ and a measure $\mu$ as follows.
If the foliation $\SF$ is defined via local charts
$\phi_i: U_i \lra \bbr^2$ (where $\{U_i\}$ is a covering
of $\SR -\{p_1, \dots, p_l\}$) which send the 
leaves of $\SF$ to horizontal arcs in $\bbr^2$, then the
transition functions 
$\phi_{ij}: \phi_i(U_i) \lra \phi_j(U_j)$
on $\phi_i(U_i)\subset \bbr^2$ 
are of the form $\phi_{ij}(x,y) = (h(x,y), c \pm y)$;
here the function $h$ is an arbitrary continuous map,
but $c$ is a constant. We require that the foliation
in a neighborhood (in $\overline{\SR}$)
of the singularities be topologically
equivalent to those that occur at the origin in $\bbc$ of the
integral curves of the line field $z^kdz^2 > 0$ where $k \ge -1$.
(There are easy extensions to arbitrary integral $k$, but we
will not need those here.) 

We define the measure $\mu$ on arcs $A\subset\SR$ as follows:
the measure $\mu(A)$ is given by
$$
\mu(A) = \int_A |dY|
$$
where $|dY|$ is defined, locally, to be the
pullback $|dY|_{U_i}= \phi_i^*(|dy|)$ of the horizontal
transverse measure $|dy|$ on $\bbr^2$.  Because of the
form of the transition functions $\phi_{ij}$ above,
this measure is then well-defined on arcs in $\SR$.

An important feature
of this measure (that follows from its definition above)
is its ``translation invariance''. That is, suppose
$A_0\subset\SR$ is an arc transverse to the foliation $\SF$,
with $\p A_0$ a pair of points, one on the leaf $l$ and one on the
leaf $l'$; then, if
we deform $A_0$ to $A_1$ via an isotopy through arcs $A_t$
that maintains the
transversality of the image of $A_0$ at every time, and also
keeps the endpoints of the arcs $A_t$ fixed on the leaves
$l$ and $l'$, respectively, then we observe that
$\mu(A_0)=\mu(A_1)$.

%Two measured foliations $(\SF,\mu)$ and $(\SG,\nu)$ are said to be
%equivalent if after possibly some Whitehead moves on $\SF$ and
%$\SG$, there is a self-homeomorphism of $\SR$ which takes $\SF$ to
%$\SG$, and $\mu$ to $\nu$. Here a Whitehead move is the
%transformation of one foliation to another by collapsing a finite
%arc of a leaf between two singularities, or the inverse procedure
%(see \cite{FLP}).

Now a holomorphic quadratic differential $\Phi$ defines a measured
foliation in the following way. The zeros $\Phi^{-1}(0)$ of $\Phi$
are well-defined; away from these zeros, we can choose a canonical
conformal coordinate $\z(z)=\int^z\sqrt\Phi$ so that $\Phi=d\z^2$.
The local measured foliations ($\{\ree\z=\oper{const}\}$,
$|d\ree\z|$) then piece together to form a measured foliation known
as the vertical measured foliation of $\Phi$, with the
translation invariance of this measured foliation of $\Phi$
following from Cauchy's theorem.

Work of Hubbard and Masur (\cite{HM}) (see also alternate proofs
in [Ke], [Gar] and [Wo]), following Jenkins (\cite{J})
and Strebel (\cite{Str}), showed that given a measured foliation
$(\SF,\mu)$ and a Riemann surface $\SR$, there is a unique holomorphic
quadratic differential $\Phi_\mu$ on $\SR$ so that the horizontal
measured foliation of $\Phi_\mu$ is equivalent to  $(\SF,\mu)$.

Extremal length. The extremal length $\ext_\SR([\g])$ of a class of
arcs $\G$ on a Riemann surface $\SR$ is defined to be the conformal
invariant
$$
\sup_\rho\f{\ell^2_\rho(\G)}{\text{Area}(\rho)}
$$
where $\rho$ ranges over all conformal metrics on $\SR$ with areas
$0<\text{Area}(\rho)<\infty$ and $\ell_\rho(\G)$ denotes the infimum
of $\rho$-lengths of curves $\g\in\G$. Here $\G$ may consist of all
curves freely homotopic to a given curve, a union of free homotopy
classes, a family of arcs with endpoints in a pair of given
boundaries, or even a more general class.  Kerckhoff (\cite{K})
showed that this definition of extremal lengths of curves extended
naturally to a defintion a extremal lengths of measured foliations.

For a class $\G$ consisting of all curves freely homotopic to a
single curve $\g\subset M$, (or more generally, a measured foliation
$(\SF, \mu)$ we see that $\ext_{(\cd)}(\G)$ (or $\ext_{(\cd)}(\mu)$)
can be construed as a real-valued function $\ext_{(\cd)}(\G)$:
Teich$(M)\lra\SR$. Gardiner (\cite{Gar}) showed that
$\ext_{(\cd)}(\mu)$ is differentiable and Gardiner and Masur ([GM])
showed
that $\ext_{(\cd)}(\mu)\in C^1$ (Teich$(M)$).
[In our particular applications, the extremal length
functions on our moduli spaces will be real analytic: this
will be explained in \S 4.5.]
Moreover Gardiner
computed that
$$
d\ext_{(\cd)}(\mu)\bigm|_{[\SR]} = 2\Phi_{\mu}
$$
so that
$$
\(d\ext_{(\cd)}(\mu)\bigm|_{[\SR]}\)[\nu] =
4\ree\int_\SR\Phi_\mu\nu.\tag2.1
$$

\tec maps, \tec distance.
(This material will only be used in an extended
digression in \S5.5) Recall that points in \tec space
can also be defined to be
equivalence classes of Riemann surface structures $R$ on $M$, the
structure $\SR_1$ being equivalent to the structure $\SR_2$ if
there is a homeomorphism $h:M\to M$, homotopic to the identity,
which is a conformal map of the structures $\SR_1$ and $\SR_2$.

We define the \tec distance $d(\{\SR_1\},\{\SR_2\})$ by
$$
d_{Teich}(\{\SR_1\}, \{\SR_2\}) = \f12\log\inf_h K(h)
$$
where $h:\SR_1\to \SR_2$ is a quasiconformal homeomorphism homotopic
to the identity on $M$ and $K[h]$ is the maximal dilatation of
$h$.
This metric is well-defined, so we may unambiguously write
$\SR_1$ for $\{\SR_1\}$.

An extraordinary fact about this metric is that the extremal
maps,
known as \tec maps, admit an explicit description, as does the
family of maps which describe a geodesic.

%Specifically, let $q\in\qd(\SR)$ denote a holomorphic quadratic
%differential on $\SR$. A horizontal (resp\. vertical) trajectory is
%an arc along which $q(z)dz^2>0$ (resp\. $q(z)dz^2<0$) except at
%the
%zeros of $q$. A trajectory is critical if it passes through a
%critical point; otherwise it is regular. If $z$ is a local
%parameter near $p\in S$ with $q(p)\neq0$ and $z(p)=z_0$, then
%$w=\int^z_{z_0}q(z)^{1/2}dz$ is the natural parameter $q$ near
%$p$.
%The line element $|q(z)|^{1/2}|dz|$ defines the $q$-metric on
%$\SR$.

\tecs theorem asserts that if $\SR_1$ and $\SR_2$ are distinct
points in $T_g$, then there is a unique quasiconformal
$h : \SR_1\to \SR_2$ with $h$ homotopic to the identity on $M$ which
minimizes
the maximal dilatation of all such $h$. The complex dilatation of
$h$ may be written $\mu(h)=k\f{\bar q}{|q|}$ for some non-trivial
$q\in\qd(\SR_1)$ and some $k$, $0<k<1$, and then
$$
d_{Teich}(\SR_1,\SR_2) = \f12\log(1+k)/(1-k).
$$
Conversely, for each $-1<k<1$ and non-zero $q\in\qd(S_1)$, the
quasiconformal homeomorphism $h_k$ of $\SR_1$ onto $h_k(\SR_2)$,
which
has complex dilatation $k\bar q/|q|$, is extremal in its homotopy
class. Each extremal $h_k$ induces a quadratic differential
$q'_k$
on $h_k(\SR_1)$, with critical points of $q$ and $q'_k$
corresponding
under $h_k$; furthermore, to the natural parameter $w$ for $q$
near $p\in S_1$ there is a natural parameter $w'_k$ near $h_k(p)$
so that
$$
\Ree w'_k = K^{1/2}\Ree w\quad\text{and}\quad
\im w'_k = K^{-1/2}\im w,
$$
where $K=(1+k)/(1-k)$.
In particular, the horizontal (and vertical) foliations for $q$
and $q'_k$ are equivalent.

The map $h_k$ is called the \tec extremal map determined by $q$
and
$k$; the differential  $q$ is called the initial differential and
the differential $q_k$ is called the
terminal differential.   We can assume all quadratic
differentials are normalized in the sense that $$||q||=\int
|q|=1.$$ The \tec geodesic segment between $S_1$ and $S_2$
consists of all points $h_s(\SR_1)$ where the $h_s$ are
\tec maps on $\SR_1$ determined by the quadratic differential
$q\in\qd(\SR_1)$ corresponding to the \tec map $h : \SR_1\to \SR_2$
and $s\in[0, \|\mu(h)\|_\infty]$.

Kerckhoff \cite{K} has given a characterization of the \tec
metric
$d_{Teich}(\SR_1, \SR_2)$ in terms of the extremal lengths of
corresponding curves on the surfaces. He proves
$$
d_{Teich}(\SR_1, \SR_2) = \f12\log\sup\lm_\g\
\f{\ext_{\SR_1}(\g)}{\ext_{\SR_2}(\g)}\tag2.2
$$
where the supremum ranges over all simple closed curves on $M$.

\newpage
\sub{\S3. From Zigzags to minimal surfaces}

Let $Z$ be a zigzag of genus $p$ dividing the plane into two regions
$\nez\Om$ and $\sww\Om$. We denote the vertices of $\nez\Om$
consecutively by $P_{-p},\ldots,P_p$ and set $P_{\infty}=\infty$. The
vertices of $\sww\Om$ however are labeled in the opposite order
$Q_j := P_{-j}$ and $Q_{\infty} =\infty$. We double both regions to
obtain punctured spheres $\nez S$ and $\sww S$ whose punctures are
also called $P_j$ and $Q_j$. Finally we take hyperelliptic covers
$\nez\SR$ over $\nez S$, branched over the $P_j$, and $\sww\SR$
over $\sww S$, branched over the $Q_j$, to obtain two hyperelliptic
Riemann surfaces of genus $p$, punctured at the \wei points which
will still be called $P_j$ and $Q_j$. The degree $2$ maps to the
sphere are called $\nez\pi:\nez\SR\to \nez S$ and
$\sww\pi:\sww\SR\to\sww S$.

\example{Example 3.1} For a genus $1$ zigzag, the Riemann surfaces
$\nez\SR$ and $\sww\SR$ will be square tori punctured at the three 
half-period points and the one full-period point.
\endexample

Now suppose that the zigzag $Z$ is reflexive. Then there is a
conformal map $\phi:\nez\Om\to\sww\Om$ such that $\phi(P_j) =
Q_j$. Clearly $\phi$ lifts to conformal maps $\phi:\nez S\to
\sww S$ and $\phi:\nez\SR\to\sww\SR$ which again take punctures to
punctures.

The surface $\nez\SR$ will be the Riemann surface on which we are
going to define the \wei data. The idea is roughly as follows: If we
look at the \wei data for the Enneper surface, it is evident that
the $1$-forms $g dh$ and $\frac1g dh$ have simpler divisors than
their linear combinations which actually appear in the \wei
representation as the first two coordinate differentials. Thus we
are hunting for these two $1$-forms, and we want to define them by
the geometric properties of the (singular) flat metrics on the
surface for which they specify the line elements, because this will
encode the information we need to solve the period problem.

To do this, we look at the flat metrics on $\SR_{\nes}$ and $\SR_{\sw}$
which come
from the following
construction. First, the domains  $\nez\Om$ and $\sww\Om$ obviously carry
the flat euclidean metric ($ds=|dz|$) metrics.
Doubling these regions defines flat (singular)
metrics on the spheres $\nez S$ and $\sww S$ (with cone points
at the lifts of the vertices of the zigzags and cone angles
of alternately $\pi$ and $3\pi/2$). These metrics are then lifted to
$\SR_{\nes}$ and $\SR_{\sw}$ by the respective covering projections.

The exterior derivatives of the multivalued
developing
maps define  single valued holomorphic nonvanishing $1$-forms $\om_{\nes}$ on
$\SR_{\nes}$ and
$\om_{\sw}$ on
$\SR_{\sw}$, because the flat metrics on the punctured surfaces
have trivial
linear holonomy. Furthermore, the behavior of these $1$-forms at a puncture is
completely determined by the cone angle of the flat metric at the puncture.
Indeed, in a suitable local coordinate, the developing map of the flat metric
near a puncture with cone angle $2\pi k$ is given by $z^k$. Hence the exterior
derivative of the developing map will have a zero (or pole) of order $k-1$
there. Note that these
considerations are valid for the point $P_{\infty}$ as well if we allow
negative cone angles. All this is well known in the context of meromorphic
quadratic differentials, see \cite{Str}.

\example{Examples 3.2} For the genus $0$ zigzag,
we obtain a $1$-form $\nez\om$ on the
sphere $\nez\SR$ with a pole of order $2$ at $P_\infty$ which we
can call $dz$ and a $1$-form $\sww\om$ on the sphere $\sww\SR$ with a
zero of order $2$ at $P_0$ and a pole of order $4$  at $P_\infty$
which then is $z^2 dz$.

For a symmetric genus $1$ zigzag, $\nez\om$ will be closely related
to the \wei $\wp$-function on the square torus, as it is a $1$-form with
double zero in $P_0$ and double pole at $P_\infty$. Furthermore,
$\sww\om$ is a meromorphic $1$-form with double order zeroes in
$P_{\pm1}$ and fourth order pole at $P_\infty$ which also can be
written down in terms of classical elliptic functions.
\endexample

In general, we can write
down the divisors of our meromorphic $1$-forms as:
$$\align
(\om_{\nes}) &= \left\{
\alignedat2
&P^2_0\cdot P^2_{\pm2}\cdot P^2_{\pm4}\cdot\ldots\cdot P^2_{\pm (p-1)} \cdot
P^{-2}_\infty &&\qquad
\hbox{p odd}
\\ &P^2_{\pm1}\cdot P^2_{\pm3}\cdot\ldots\cdot P^2_{\pm (p-1)} \cdot
P^{-2}_\infty
&&\qquad
\hbox{p even}
\endalignedat \right. \\
(\om_{\sw}) &= \left\{
\alignedat2
&Q^2_{\pm1}\cdot Q^2_{\pm3}\cdot\ldots\cdot Q^2_{\pm p} \cdot Q^{-4}_\infty
&&\qquad
\hbox{p odd}
\\ & Q^2_0\cdot Q^2_{\pm2}\cdot Q^2_{\pm4}\cdot\ldots\cdot Q^2_{\pm p} \cdot
Q^{-4}_\infty &&\qquad
\hbox{p even}
\endalignedat \right. \\
\vspace{2\jot
}(d\pi_{\nes}) &=P_0\cdot P_{\pm1}\cdot P_{\pm2}\cdot\ldots\cdot
P_{\pm p} \cdot P^{-3}_\infty
\endalign$$

Now denote by $\alpha := e^{-\pi i/4}\cdot\om_{\nes},\  \beta := e^{-\pi
i/4}\cdot
\phi^*\om_{\sw}$ and
$dh := const\cdot d\pi_{\nes}$,
where we choose the constant such that
$\alpha\cdot\beta=dh^2$ which is possible because the divisors coincide. Now we
can write down the Gau{\ss} map of the \wei data on $\SR_{\nes}$ as $g=
\frac\alpha{dh}$, and we check easily that the line element (2.0)
is regular everywhere on $\SR_{\nes}$ except at the lift of $P_\infty$.

One can check that the thus defined \wei data coincide for the reflexive
genus $0$ and genus $1$ zigzags with the data for the Enneper surface and the
Chen-Gackstatter surface.

We can now claim

\sub{Theorem 3.3} If $Z$ is a symmetric reflexive zigzag of genus $p$, then
$(\SR_{\nes}, g, dh)$ as above define a \wei representation of a minimal
surface
of genus $p$ with one Enneper-type end and total curvature $-4\pi(p+1)$.

\sub{Proof}
The claim now is that the $1$-forms in the \wei representation all
have purely imaginary periods. For $dh$ this is obvious,
because the form $dh$ is even exact and so all periods even vanish.
Because of
$(g-\frac1g)dh = \alpha-\beta$ and
$i(g+\frac1g)dh = i(\alpha+\beta)$ this is equivalent to the claim
that $\alpha$ and $\beta$ have complex conjugate periods. To see
this, we first construct a basis for the homology on $\nez\SR$ and
then compute the periods of $\alpha$ and $\beta$ using their
geometric definitions.

To define $2p$ cycles $B_j$ on $\nez\SR$, we take curves
$b_j$ in $\nez\Om$
connecting a boundary point slightly to the right of
$P_{j+1}$ with a boundary point slightly to the left of $P_j$  for
$j=-p,\ldots, p-1$. We double this curve to obtain a closed curve
$B_j$ on $\nez\SS$ which encircles exactly $P_j$ and
$P_{j+1}$. These curves have closed
lifts $B_j$ to $\nez\SR$ and form a homology basis. Now to
compute a period of our $1$-forms, observe that a period is nothing
other than the image of the closed
curve under the developing map of the flat metric which defines the
$1$-form. But this developing map can be read off from the zigzag
--- one only has to observe that developing a curve around a vertex
(regardless whether the angle there is $\frac\pi2$ or $3\frac\pi2$)
will change the direction of the curve there by $180^\circ$. Doing
this yields
$$
\align
\int_{B_j}\alpha &= \int_{B_j} e^{-\pi i/4}\cdot\nez\om
= 2 e^{-\pi i/4}\cdot (P_j-P_{j+1})
\\
\int_{B_j} \beta &= \int_{B_j} e^{-\pi i/4}\cdot
\phi^*\sww\om =
\int_{\phi(B_j)}e^{-\pi i/4}\cdot\sww\om
= 2 e^{-\pi i/4}\cdot (Q_j-Q_{j+1}) \\
&= 2 e^{-\pi i/4}\cdot (P_{-j}-P_{-j-1})
\endalign
$$
which yields the claim by the symmetry of the zigzag.

Finally, we have to compute the total absolute curvature of the minimal
surface.
By the definition of the Gau{\ss} map we have

$$
(g) = \left\{
\alignedat2
&P_0^{+1}\cdot P^{-1}_{\pm1}\cdot P^{+1}_{\pm2}\cdot\ldots\cdot P^{-1}_{\pm
p} \cdot P_\infty
&&\qquad
\hbox{p odd}
\\ &P^{-1}_0\cdot P^{+1}_{\pm1}\cdot P^{-1}_{\pm2}\cdot\ldots\cdot
P^{-1}_{\pm p} \cdot P_\infty
&&\qquad
\hbox{p even}
\endalignedat \right.
$$
and thus $\oper{deg} g = p+1$ which implies $\int_R K dA = -4\pi (p+1)$
as claimed.
\qed

\sub{Remark 3.4}
We close by making some comments on the amount of symmetry involved
in this approach. Usually in the construction of minimal surfaces
the underlying Riemann surface is assumed to have so many
automorphisms that the moduli space of possible conformal structures
is very low dimensional (in fact, it consists only of one point in
many examples). This helps solving the period problem (if it is
solvable) because this will then be a problem on a low dimensional
space. In our approach, the dimension of the moduli space grows with
the genus, and the use of symmetries has other purposes: It allows
us to construct for given periods a pair of surfaces with one
$1$-form on each which would solve the period problem if only the
surfaces would coincide.
Indeed we observe

\proclaim{Lemma 3.5} The minimal surface
of genus $p$ constructed below has only eight
isometries, and at most eight conformal or anticonformal
automorphisms that fix the end, independently of genus.
\endproclaim

\sub{Proof} Observe that as the end is unique, any isometry of the
minimal surface fixes the end. As the isometry necessarily
induces a conformal or anti-conformal automorphism, it is
sufficient to prove only the latter statement of the lemma.
Because of the uniqueness of the hyperelliptic involution, this
automorphism descends to an automorphism of the punctured sphere
which fixes the image of infinity and permutes the punctures
(the images of the \wei points). As there are at least three
punctures, all lying on the real line, we see that the real line
is also fixed (setwise).  After taking the reflection
(an anti-conformal automorphism) of the
sphere across the real line, which fixes all the punctures, we are
left to consider the conformal automorphisms of the domain $\nez\Om$.
Finally, this domain
has only two conformal symmetries, the identity and the reflection about the
diagonal.  The lemma follows by counting the automorphisms we have
identified in the discussion. \qed

\newpage
\sub{\S4. The Height Function on Moduli Space}

\sub{4.1 The space of zigzags $\SZ_p$ and a natural compactification
$\p\ov\SZ$} We recall the space $\SZ_p$ of
equivalence classes of symmetric
genus~$p$ zigzags constructed in Section~2.2; here the equivalence
by similarity was defined so that two zigzags $Z$ and $Z'$ would be
equivalent if and only if both of the pairs of complementary
domains $(\nez\Om(Z),\nez\Om(Z'))$ and
$(\sww\Om(Z),\sww\Om(Z'))$ were conformally equivalent. Label the
finite vertices of the zigzag by $P_{-p},\dots,P_0,\dots,P_p$. Thus,
we may choose a unique representative for each class in $\SZ_p$ by
setting the vertices $P_0\in\{y=x\}$,  $P_p=1$, $P_{-p}=i$ and
$P_k=i\ov{P_{-k}}$ for $k=0,\dots,p$; here all the vertices
$P_{-p},\dots,P_p$ are required to be distinct. The topology of
$\SZ_p$ defined in Section~2.2 then agrees with the topology of the
space of canonical representatives induced by the embedding of
$\SZ_p\to\BC^{p-1}$ by $Z\mapsto(P_1,\dots,P_{p-1})$. With these
normalizations and this last remark on topology, it is evident that
$\SZ_p$ is a cell of dimension $p-1$.

We have interest in the natural compactification of this cell,
obtained by attaching a boundary $\p\ov\SZ_p$ to $\SZ$. This
boundary will be composed of zigzags where some proper consecutive
subsets of $\{P_0,\dots,P_p\}$ (and of course the reflections of
these subsets across $\{y=x\}$) are allowed to coincide; the
topology on $\ov\SZ_p=\SZ_p\cup\p\ov\SZ_p$ is again given by the
topology on the map of coordinates of normalized representatives
$Z\in\SZ_p\mapsto(P_1,\dots,P_{p-1})\in C^{p-1}$.

Evidently, $\p\ov\SZ_p$ is stratified by unions of zigzag spaces
$\SZ^k_p$ of real dimension $k$, with each component of $\SZ^k_p$
representing the (degenerate) zigzags that result from allowing $k$
distinct vertices to remain in the (degenerate) zigzag after some
points $P_0,\dots,P_p$ have coalesced. For instance $\SZ^0_p$
consists of the zigzags where all of the points $P_0,\dots,P_p$ have
coalesced to either $P_0\in\{y=x\}$ or $P_1=1$, and the faces $\SZ^{p-2}_p$
are the loci in $\ov\SZ_p$ where only two consecutive vertices have
coalesced.

Observe that each of these strata is naturally a zigzag space in
its own right, and one can look for a reflexive symmetric zigzag of
genus~$k+1$ within $\SZ^k_p$.

\bpage
\sub{4.2 Extremal length functions on $\SZ$} Consider the
punctured sphere $S_{\nes}$ in \S3, where we labelled the punctures
$P_{-p},\dots,P_0,\dots,P_p$, and $P_\infty$ and observed that
$\nez S$ had two reflective symmetries: one about the image of
$Z$ and one about the image of the curve $\{y=x\}$ on $\nez\Om$.
Let $[B_k]$ denote the homotopy class of simple curves which
encloses the punctures $P_k$ and $P_{k+1}$ for $k=1,\dots,p-1$ and
$[B_{-k}]$ the homotopy class of simple curves which encloses the
punctures $P_{-k}$ and $P_{-k-1}$ for $k=1,\dots,p-1$. Let $[\g_k]$
denote the pair of classes $[B_k]\cup[B_{-k}]$. Under the homotopy
class of maps which connects $\nez S$ to $\sww S$ (lifted from
$\phi:\nez\Om\to\sww\Om$, the vertex preserving  map), there are
corresponding homotopy classes of curves on $\sww S$, which we will
also label $[\g_k]$.

Set $\nez E(k)=\ext_{\nez S}([\g_k])$ and
$\sww E(k)=\ext_{\sww S}([\g_k])$ denote the extremal lengths of
$[\g_k]$ in $\nez S$ and $\sww S$, respectively.

Let $\tsym T$ denote a subspace of the \tec space of $2p+2$ punctured
spheres whose points are equivalence classes of $2p+2$ punctured
spheres (with a pair of involutions) coming from one complementary
domain $\nez\Om(Z)$ of a symmetric zigzag $Z$. This $\tsym T$ is
a $p-1$ dimensional subspace of the \tec space $\tp$ of $2p+2$
punctured spheres.

Consider the map $\nez E:\tsym T\to\BR^{p-1}_+$ given by
$\nez S\mapsto(\nez E(1),\dots,\nez E(p-1))$.

\bpage
\proclaim{Proposition 4.2.1} The map $\nez E:\tsym T\to\BR^{p-1}_+$
is a homeomorphism onto $\BR^{p-1}_+$.
\endproclaim

\sub{Proof} It is clear that $\nez E$ is continuous. To see
injectivity and surjectivity, apply a Schwarz-Christoffel map
$SC: \nez\Om\to\{\imm z>0\}$ to $\nez\Om$; this map sends $\nez\Om$
to the upper half-plane, taking $Z$ to $\BR$ so that
$SC(P_\infty)=\infty$, $SC(P_0)=0$, $SC(P_1)=1$ and
$SC(P_{-k})=-SC(P_k)$. These conditions uniquely determine $SC$;
moreover $\nez E(k)=2\ext_{\BH^2}(\G_k)$ where $\G_k$ is the class
of pairs of curves in $\BH^2$ that connect the real arc between
$SC(P_{-k-2})$ and $SC(P_{-k-1})$ to the real arc between
$SC(P_{-k})$ and $SC(P_{-k+1})$, and the arc between $SC(P_{k-1})$
and $SC(P_k)$ to the arc between $SC(P_{k+1})$ and $SC(P_{k+2})$.
Now, any choice of $p-1$ numbers $x_i=SC(P_i)$ for $2\le1\le p$
uniquely determines a point in $\tsym T$, and these choices are
parametrized by the extremal lengths
$\ext_{\BH^2}(\G_k)\in(0,\infty)$. This proves the result. \qed

Let $\qd^{\text{symm}}(\nez S)$ denote the vector space of
holomorphic quadratic differentials on $\nez S$ which have at
worst simply poles at the punctures and are real along the image of
$Z$ and the line $\{y=x\}$.

Our principal application of Proposition 4.2.1 is the following

\bpage
\proclaim{Corollary 4.2.2} The cotangent vectors $\{d\nez E(k)
\mid k=1,\dots,p-1\}$ (and $\{d\sww E(k)\mid k=1,\dots,p-1\}$)
are a basis for $T^*_{\nez S}\tsym T$, hence for
$\qd^{\text{symm}}(\nez S)$.
\endproclaim

\sub{Proof} The cotangent space $T^*_{\nez S}\tsym T$ to the \tec
space $\tsym T$ is the space $\qd(\nez S)$ of holomorphic quadratic
differentials on $\nez S$ with at most simple poles at the
punctures. A covector cotangent to $\tsym T$ must respect the
reflective symmetries of the elements of $\tsym T$, hence its
horizontal and vertical foliations must be either parallel or
perpendicular to the fixed sets of the reflections. Thus such a
covector, as a holomorphic quadratic differential, must be real on
those fixed sets, and hence must lie in
$\qd^{\text{symm}}(\nez S)$. The result follows from the
functions $\{\nez E(k)\mid k=1,\dots,p-1\}$ being coordinates for
$\tsym T$. \qed

\bpage
\sub{4.3  The height function $D(Z):\SZ\to\BR$}
Let the height function $D(Z)$ be
$$
D(Z) = \sum^{p-1}_{j=1}\[\exp\(\f1{E_{\nes}(j)}\) -
\exp\(\f1{E_{\sw}(j)}\)\]^2 +\[E_{\nes}(j) - E_{\sw}(j)\]^2.\tag4.1
$$
We observe that $D(Z)=0$ if and only if $E_{\nes}(j)=E_{\sw}(j)$,
which holds if and only if $\nez S$ is conformally equivalent to
$\sww S$. We also observe that, for instance, if
$E_{\nes}(j)/E_{\sw}(j)\ge C_0$ but both $E_{\nes}(j)$ and
$E_{\sw}(j)$ are quite small, then $D(Z)$ is quite large. It is this
latter fact which we will exploit in this section.

\sub{4.4 Monodromy Properties of the Schwarz-Christoffel Maps}
Here we derive the facts about the Schwarz-Christoffel maps we need
to prove properness of the height function.

Let $t=(0<t_1<t_2< \ldots < t_p)$ be $p$ points on the real line. We
put $t_0:=0, t_\infty:=\infty$ and $t_{-k}=-t_k$. Then the
Schwarz-Christoffel formula tells us that we can map the upper half
plane conformally to a NE-domain such that the $\{t_i\}$ are mapped to
vertices by the function

$$
f(z) = \int_0^z (t-t_{-p})^{1/2}(t-t_{-p+1})^{-1/2}\cdots
(t-t_p)^{1/2} dt
$$
and to a SW-domain by
$$
g(z) = \int_0^z (t-t_{-p})^{-1/2}(t-t_{-p+1})^{+1/2}\cdots
(t-t_p)^{-1/2} dt
$$

Note that the exponents alternate sign. We are not interested in
normalizing these maps at the moment by introducing some factor, but
we have to be aware of the fact that scaling the $t_k$ will scale
$f$ and $g$.

Now introduce the periods $a_k = f(t_{k+1})-f(t_k)$ and
$b_k=g(t_{k+1})-g(t_k)$ which are complex numbers, either real or
purely imaginary. Denote by
$$
\Ts:=\{t:t_i\in\Bbb C, t_j\ne t_k\quad \forall j,k\},
$$
the complex-valued configuration space for the $2p+1$-tuples $\{t\}$.
It is clear that we can analytically continue the $a_k$ and $b_k$
along any path in $\Ts$ to obtain holomorphic multi-valued functions.

\proclaim{Lemma 4.4.1} Continue $a_k$ analytically along a path in
$T$ defined by moving $t_j$ anticlockwise around $t_{j+1}$ and denote
the continued function by $\tilde a_k$, similarly for $b_k$. Then we
have
$$
\tilde a_k = \left\{\alignedat2
  &a_k &&\quad\hbox{if}\  k\ne j-1,j+1 \\
  &a_k+2a_{k+1} &&\quad\hbox{if}\ j=k+1\\
  &a_k-2a_{k-1} &&\quad\hbox{if}\ j=k-1\\
\endalignedat \right.
$$
with analogous formulas holding for $\tilde b_k$.
\endproclaim

\sub{Proof} Imagine that the defining paths of integration for
$a_k$ was made of flexible rubber band which is tied to  $t_k$ $t_{k+1}$.
Now moving $t_j$ will possibly drag the rubber band into some new position. The
resulting curves are precisely those paths of integration which need to be used
to compute $\tilde a_k$. If $j \ne k-1,k+1$, the paths remain the same,
hence $\tilde a_k=a_k$. If $j=k+1$, the rubberband between $t_k$ and
$t_{k+1}$ is pulled around $t_{k+2}$ and back to $t_{k+1}$. Hence
$a_k$ changes by the amount of the integral which goes from
$t_{k+1}$ to $t_{k+2}$, loops around $t_{k+2}$ and then back to
$t_{k+1}$. Hence the first part contributes $a_{k+1}$. Now,
for the second part of the path of integration,  by the
very definition of the Schwarz-Christoffel maps we know that a small
interval through
$t_{k+2}$ is mapped to a $90^\circ$ hinge, so that a
small infinitesimal loop turning around $t_{k+2}$ will be mapped to
an infinitesimal straight line segment. In fact, locally near $t_{k+2}$ the
Schwarz-Christoffel map is of the form $z \mapsto z^{1/2}$ or  $z \mapsto
z^{3/2}$.
Therefore we get from the integration back to $t_{k+1}$ another contribution of
$+a_{k+1}$. The same argument is valid for $j=k-1$ and for the $b_k$. \qed

%\vskip .5cm
%\nointerlineskip
%\moveright0cm\hbox{\vrule height  1cm $P_{k-1}$ \hskip 1cm $P_k$}
%\nointerlineskip
%\moveright0cm \vbox{\hrule width 2cm}
%\nointerlineskip
%\moveright2cm \hbox{\vrule height 0.5cm $P_{k+1}$ \hskip 1.5cm $P_{k+2}$}
%\nointerlineskip
%\moveright2cm \vbox{ \hrule width 2.5cm }
%\nointerlineskip
%\moveright4.5cm \hbox{\vrule height 1cm}
%\nointerlineskip

Now denote by $\delta:= t_{k+1}-t_k  $ and fix all $t_j$ other than
$t_{k+1}$: we regard $t_{k+1}$ as the independent variable.

\proclaim{Corollary 4.4.2} The functions
$a_k-\frac{\log\delta}{\pi i} a_{k+1}$ and
$b_k-\frac{\log\delta}{\pi i} b_{k+1}$ are holomorphic in $\delta$ at
$\delta=0$.
\endproclaim

\sub{Proof}  By Lemma 4.4.1, the above functions are singlevalued
and holomorphic in a punctured neighborhood of $\delta=0$. It is
easy to see from the explicit integrals
defining $a_k$, $b_k$, $a_{k+1}$ and $b_{k+1}$
that the above functions are
also bounded, hence they extend holomorphically to $\delta=0$. \qed

Now for the properness argument, we are more interested
in the absolute values of the periods than in the periods
themselves:
we translate the above statement about periods into a statement
about their respective absolute values.
This will lead to a crucial difference in the behavior
 of the extremal length functions on
the NE-- and SW regions.

\proclaim{Corollary 4.4.3} Either $\vert a_k\vert - \frac{\log\delta}\pi \vert
a_{k+1} \vert$ or $\vert a_k\vert + \frac{\log\delta}\pi \vert
a_{k+1} \vert$ is real analytic in $\delta$ for $\delta=0$. In the first case,
$\vert b_k\vert + \frac{\log\delta}\pi \vert
b_{k+1} \vert$ is real analytic in $\delta$, in the second
$\vert b_k\vert - \frac{\log\delta}\pi \vert
b_{k+1} \vert$.
\endproclaim

\sub{Remark} Note the different signs here! This reflects that we
alternate between left and right turns in the zigzag.

\sub{Proof} If we follow the images of the $t_k$ in the NE-domain,
we turn alternatingly left and right, that is, the
direction of $a_{k+1}$ alternates between $i$ times the direction of
$a_k$ and $-i$ times the direction of
$a_k$.

This proves the
first statement, using corollary 4.4.2. Now if we turn left at $P_k$
in the NE domain, we turn right in the SW domain, and vice versa,
because the zigzag is run through in the opposite orientation. This
proves the second statement. \qed

From this we deduce a certain non-analyticity which is used in the properness
proof. Denote by $s_k, t_k$ the preimages of the vertices $P_k$ of a zigzag
under the Schwarz-Christoffel maps for the NW- and SW-domain respectively. We
normalize these maps such that
$s_0=t_0=0$ and $s_n=t_n=1$. Introduce $\delta_{\nes}=s_{k+1}-s_k$ and
$\delta_{\sw}=t_{k+1}-t_k$. We can now consider $\delta_{\nes}$ as a function
of $\delta_{\sw}$:

\proclaim{Corollary 4.4.4} The function $\delta_{\nes}$ does not depend
real analytically on $\delta_{\sw}$.
\endproclaim
\demo{Proof}
%Suppose the opposite. Then by corollary 4.4.3, either
%$\vert b_k\vert + \frac{\log\delta_{\sw}}\pi \vert
%b_{k+1} \vert$ or $\vert b_k\vert - \frac{\log\delta_{\sw}}\pi \vert
%b_{k+1} \vert$  depends real analytically on $\delta_{\sw}$. By the
%definition of $\delta_{\sw}$ and $\delta_{\nes}$ and again corollary
%4.4.3 in fact both would depend real analytically on $\delta_{\sw}$
%so that also
%$\log\delta\vert b_{k+1}\vert$ would depend real analytically on
%$\delta_{\sw}$. But by the explicit Schwarz-Christoffel integrals, we
%know that $\vert b_{k+1}\vert$
%does depend real analytically on $\delta_{\sw}$ which is a
%contradiction.
Suppose the opposite is true.We know that either
$\vert b_k\vert + \frac{\log\delta_{\sw}}\pi \vert
b_{k+1} \vert$ or $\vert b_k\vert - \frac{\log\delta_{\sw}}\pi \vert
b_{k+1} \vert$  depends real analytically on $\delta_{\sw}$,
hence on $\delta_{\nes}$, so we may assume with no loss in
generality that $\vert b_k\vert + \frac{\log\delta_{\sw}}\pi \vert
b_{k+1} \vert$ depends real analytically on $\delta_{\sw}$. Then
by Corollary 4.4.3, we see that
$\vert a_k\vert - \frac{\log\delta_{\sw}}\pi \vert
a_{k+1} \vert$ depends analytically on $\delta_{\nes}$, hence by
assumption on $\delta_{\sw}$.  Hence $B(\delta_{\sw}) :=
\frac{\vert b_k\vert}{\vert b_{k+1}\vert} +
\frac{\log\delta_{\sw}}\pi$
and $A(\delta_{\nes}):=\frac{\vert a_k\vert}{\vert a_{k+1}\vert} -
\frac{\log\delta_{\nes}}\pi$ depends real analytically
on $\delta_{\sw}$ and $\delta_{\nes}$, respectively.  But
$$
 \frac{\vert b_k\vert}{\vert b_{k+1}\vert} =
\frac{\vert a_k\vert}{\vert a_{k+1}\vert}
$$
by the assumption on equality of periods, so
$$
B(\delta_{\sw}) - \frac{\log\delta_{\sw}}\pi =
A(\delta_{\nes}) +\frac{\log\delta_{\nes}}\pi.
$$
But then $\frac{\log(\delta_{\sw}\delta_{\nes})}\pi
= B(\delta_{\sw}) - A(\delta_{\nes})$
is analytic in $\delta_{\nes}$.  Of course, the product
$\delta_{\sw}\delta_{\nes}$ is analytic in $\delta_{\nes}$
and non-constant, as $\delta_{\sw}\delta_{\nes}$ tends to zero by
the hypothesis on extremal length.  But then
$log(\delta_{\sw}\delta_{\nes})$ is analytic in $\delta_{\nes}$,
near $\delta_{\nes}=0$, which is absurd.
\qed
\enddemo

\sub{Remark}
Note that the Corollary remains true if we consider zigzags which
turn alternatingly left and right by a (fixed) angle
other than $\pi/2$. This will only
affect constants in Lemma 4.4.1 and Corollaries 4.4.2, 4.4.3. In Corollary
4.4.4, we need only that the coefficients of
$\log\delta$ are distinct, and this is also
true for the zigzags with non-orthogonal sides.
We will use this generalization in section 6.

\sub{4.5 An extremal length computation}
Here we compute the extremal length of curves separating two points
on the real line. This will be needed in the next section. We do
this first in a model situation: Let $\lambda<0<1$ and consider the
family of curves $\Gamma$ in the upper half plane joining the
interval $(-\infty,\lambda)$ with the interval $(0,1)$. For a
detailed account on this, see \cite{Oht}, p.\ 179--214. He gives
the result in terms of the Jacobi elliptic functions from which
it is straightforward to
deduce the asymptotic expansions which we need.
Because it fits in the spirit of this paper, we
give an informal description of what is involved in terms of
elliptic integrals of \wei type.

It turns out that the extremal metric for $\Gamma$ is rather explicit
and can be seen best by considering a slightly different problem:
Consider the family $\Gamma'$ of curves in $S^2$ encircling only
$\lambda$ and $0$ thus separating them from $1$ and $\infty$. Then
the extremal length of $\Gamma'$ is twice the one we want.

\proclaim{Lemma 4.5.1} The extremal metric in this situation is
given by the flat cone metric on $S^2-\{\lambda,0,1,\infty\}$ with
cone angles $\pi$ at each of the four vertices.
\endproclaim

\sub{Proof} Directly from the length-area method of Beurling,
or see \cite{Oht}.
\qed

This metric can be constructed by taking the double cover
over
$S^2-\{\lambda,0,1,\infty\}$, branched over
$\lambda,0,1,\infty$ which is a torus $T$, which has a unique flat
conformal metric (up to scaling). This metric descends as the cone
metric we want to $S^2-\{\lambda,0,1,\infty\}$. This allows us to
compute the extremal length in terms of certain elliptic period
integrals. Because the covering projection $p:T\to S^2$ is given by
the equation ${p'}^2 = p(p-1)(p-\lambda)$ we compute the periods of
$T$ as
$$
\omega_i = \int_{\gamma_i} \frac{du}{\sqrt {u(u-1)(u-\lambda)}}
\tag4.2
$$
where $\gamma_1$ denotes a curve in $\Gamma'$ and $\gamma_2$ a
curve circling around $\lambda$ and $0$. We conclude that the
extremal length we are looking for is given by

\proclaim{Lemma 4.5.2}
$$
\ext(\Gamma) = 2
\frac{\vert\omega_1\vert^2}{\det(\omega_1,\omega_2)}
\tag4.3
$$
\endproclaim
\sub{Proof} see \cite{Oht}
\qed

Alternatively,
$$
\frac{\omega_2}{\omega_1} = - \frac{\int_0^\lambda \frac{du}{\sqrt
{u(u-1)(u-\lambda)}}}
{\int_0^1 \frac{du}{\sqrt {u(u-1)(u-\lambda)}}}
\tag4.4
$$
This is as explicit as we can get.

It is evident from formulas (4.2), (4.3), and (4.4) that $\ext(\G)$
is real analytic on $T_{0,4}$, hence on $T^{\hbox{symm}}_{0,2p+2}$.

Now we are interested in the
asymptotic behavior of the extremal lengths $\ext(\Gamma)$ as $\lambda\to0$.

\proclaim{Lemma 4.5.3} $$\ext(\Gamma) = O\left(\frac1{\log\vert
\lambda\vert}\right)$$
\endproclaim

\sub{Proof} The asymptotic behavior of elliptic integrals is well
known, but it seems worth observing that all the information which
we need is in fact contained in the geometry. For a more formal treatment, we
refer to \cite{Oht, Rain}. The period
$\omega_1$ is easily seen to be holomorphic in $\lambda$ by developing the
integrand into a power series and integrating term by term; explicitly we can
obtain
$$
\omega_1(\lambda) = 2\pi\left(1+\frac14\lambda +
\frac9{32}\lambda^2+\frac{25}{128}\lambda^3+\ldots\right)
$$
but all we need is the holomorphy and $\omega_1(0)\ne0$.
Concerning $\omega_2$, the general theory of ordinary differential
equations with regular singular points predicts that any solution
$\omega$ of the o.d.e. has the general form
$$
\omega= c_1\big(\log(\lambda) \omega_1(\lambda) + \lambda
f_1(\lambda)\big) + c_2 \omega_1(\lambda)
$$
with some explicitly known holomorphic function $f_1(\lambda)$.

From a similar monodromy argument as in the above section 4.4 one can
obtain that
$$
\omega_2(\lambda) - \frac{i}{\pi}\log(-\lambda)\omega_1(\lambda)
$$
is holomorphic at $\lambda=0$ which is simultaneously a more
specific but less general statement. Nevertheless this can already
be used to deduce the claim, but for the sake
of completeness
we cite from \cite{Rain} the full expansion:

$$
\omega_2(\lambda) =
\frac{i}{\pi}\left(\log(-\lambda)\omega_1(\lambda)+\lambda
f_1(\lambda)\right) -\frac{4i\log2}{\pi}\omega_1(\lambda)
$$
with the expansion of $\lambda f_1(\lambda)$ given by
$$
\lambda f_1(\lambda) = \sum_{n=1}^\infty \frac{(a)_n(b)_n}{(n!)^2}
\left(H(a,n)+H(b,n)-2H(1,n)\right) \lambda^n
$$
Here
$$
\align
(a)_n &:= a(a+1)\cdot\ldots\cdot(a+n-1) \\
H(a,n) &= \frac1a+\frac1{a+1}+\cdots+\frac1{a+n-1} \\
a&=1/2\\
b&=1/2
\endalign
$$
From this, one can deduce the claim in any desired degree of
accuracy. \qed

Now we generalize this to 4 arbitary points $t_1<t_2<t_3<t_4$ on the
real line and to the family $\Gamma$ of curves connecting
the arc $t_1 t_2$ to the arc $t_3 t_4$.

We denote the cross ratio of $t_1$, $t_2$, $t_3$, $t_4$ by
$(t_1:t_2:t_3:t_4)$ which is chosen so that $(\infty: \lambda: 0:
1)=\lambda)$.

\proclaim{Corollary 4.5.4} For $t_2 \to t_3$, we have
$$
\ext \Gamma =  O\left(\frac1{-\log\vert
(t_1:t_2:t_3:t_4)\vert}\right)
$$
\endproclaim

\sub{Proof} This follows by applying the M\"obiustransformation to
the $t_i$ which maps them to $\infty, \lambda, 0, 1$. \qed

\sub{Remark} Here we can already see that to establish properness we need
to consider points $t_i$ such that $t_2\to t_3$ while $t_1$ and
$t_4$ stay at finite distance away.

\sub{4.6 Properness of the Height Function}
In this section we prove

\proclaim{Theorem 4.6.1} The height function $D(Z)$ is proper on
$\Cal Z$, for $p>2$.
\endproclaim

\sub{Proof}
Let $Z_0$ be a zigzag in the boundary of  $\ov\SZ$. We can imagine
$Z_0$ as an ordinary zigzag where some (consecutive)
vertices have coalesced. We can
assume  that we have a cluster of coalesced points
$P_k=P_{k+1}=\ldots=P_{k+l}$ but $P_{k-1}\ne P_k$ and $P_{k+l}\ne
P_{k+l+1}$
(here $k \ge 0$ and $k+l \le p$).

We first consider the case where $k\ge 1$, taking up the case
$k=0$ later.
Denote the family of curves connecting the
segment $P_{k-1}P_k$ with the segment $P_{k+l} P_{k+l+1}$
(and their counterparts symmetric about the central point $P_0$)
in the NE
domain by $\Gamma_{\nes}$ and in the SW domain by $\Gamma_{\sw}$
and their extremal lengths in $\nez\Om$ $(\sww\Om, resp.)$ by
$\ext\Gamma_{\nes}$ $(\ext\Gamma_{\sw}, resp.)$.
Now recall that the height function was defined so that
$$
D(Z) \ge \sum^{p-1}_{j=1}\[\exp\(\f1{E_{\nes}(j)}\) -
\exp\(\f1{E_{\sw}(j)}\)\]^2.
$$
where the extremal length were taken of curves encircling two consecutive
points. To prove properness it is hence sufficient to prove that at least one
pair $E_{\nes}(j), E_{\sw}(j)$ approaches $0$ with different rates to
some
order for any sequence of zigzags $Z_n\to Z_0$. Suppose this is not the case.
Then especially all $E_{\nes}(j), E_{\sw}(j)$ with $j=k,\ldots,k+l$ approach
zero at the same rate for some $Z_n\to Z_0$. Now conformally the points
$t_k,\ldots,t_{k+l+1}$ are determined by the extremal lengths
$E_{\nes}(j)$ and
$E_{\sw}(j)$ so that under the assumption, $\ext\Gamma_{\nes}$ and
$\ext\Gamma_{\sw}$ approach zero with the same rate. Thus to obtain a
contradiction it is sufficient to
prove that
$$
e^{1/\ext\Gamma_{\nes}}-e^{1/\ext\Gamma_{\sw}}
$$
is proper in a neighborhood of $Z_0$ in $\Cal Z$. Such a
neighborhood is given by
 all zigzags $Z$ where
distances between coalescing points are sufficiently
small. Especially,
the quantity $\epsilon:= \vert P_k - P_{k+j}\vert$ is small.

To estimate the extremal length, we map the NE-- and SW domains of a
zigzag $Z$ in a neighborhood of $Z_0$ by the inverse
Schwarz-Christoffel maps to the upper half plane and apply then the
asymptotic formula of the last section, using that the
asymptotics for a symmetric pair of degenerating curve families
agree with the asymptotics of a single degenerating curve family.

Denote by $\delta_{\nes}$ and $\delta_{\sw}$ the difference
$t_{k+1}-t_k$ of the images of $P_{k+1}$ and $P_k$ for NE and SW
respectively. Because the Schwarz-Christoffel map is a homeomorphism
on the compactified domains, the quantities
$\delta_{\nes}$ and $\delta_{\sw}$ will
go to zero with $\e$,
while the distances $t_{k-1}-t_k$ and
$t_{k+1}-t_{k+l+1}$ are uniformly bounded away from zero in any
compact coordinate patch. Hence by Corollary 4.5.4
$$
\align
\Big\vert
e^{\frac1{\ext\Gamma_{\nes}}}&-e^{\frac1{\ext\Gamma_{\sw}}}
\Big\vert \\ &=\Bigg\vert
O\left(\frac1{\vert(t_{k-1}:t_k:t_{k+1}:t_{k+l+1})_{\nes}\vert}\right) -
O\left(\frac1{\vert(t_{k-1}:t_k:t_{k+1}:t_{k+l+1})_{\sw}\vert}\right)\Bigg\vert
\\ &=O(\frac1{\delta_{\nes}})-O(\frac1{\delta_{\sw}})\tag4.5
\endalign
$$
On the other hand, by corollary 4.4.4, we know that $\delta_{\nes}$
cannot depend analytically on $\delta_{\sw}$ so that one term will
dominate the other and no cancellation occurs. Finally all occuring
constants are uniform in a coordinate patch and $\delta_{\nes}$
depends there in a uniform way on $\epsilon$. This proves local
properness near $Z_0$ in this case and gives the desired
contradiction.

In the case where $P_k=0, ...,P_{k+l}$ are coalescing (here $k+l <
p$), we use the other terms in the height function, i.e. the
inequality
$$
D(Z) \ge \sum^{p-1}_{j=0}\[E_{\nes}(j) - E_{\sw}(j)\]^2\tag4.6
$$
It is a straightforward exercise in the definition of extremal length
(see lemma 4.5.1) that, since $\gamma_{k+l-1}$ and $\gamma_{k+l}$
intersect once (geometrically), and the point $P_{k+l+1}$ converges
to a finite point distinct from $P_k, ...,P_{k+l}$, we can conclude
that $E_{\nes}(k+l) = O(\frac1{E_{\nes}(k+l-1)})$ and
$E_{\sw}(k+l) = O(\frac1{E_{\sw}(k+l-1)})$.
(Here we use the hypothesis that $p>2$ in the final
argument to ensure the existence of a second dual curve.)
Yet an examination
of the argument above (see also Corollary 4.4.4, especially) shows
that $E_{\nes}(k+l-1)$ vanishes at a different rate than
$E_{\sw}(k+l-1)$, hence $E_{\nes}(k+l)$ grows at a different
rate than $E_{\sw}(k+l)$. This term alone then in inequality (4.6)
shows the claim.
\qed

\newpage
\sub{\S5. The Gradient Flow for the Height Function}

\sub{5.1} To find a zigzag $Z$ for which $D(Z)=0$, we imagine
flowing along the vector field $\grad D$ on $\SZ$ to a minimum
$Z_0$. To know that this minimum $Z_0$ represents a reflexive
zigzag (i.e., a solution to our problem), we need to establish
that, at such a minimum $Z_0$, we have $D(Z)=0$. That result is the
goal of this section; in the present subsection, we state the
result and begin the proof.

\bpage
\proclaim{Proposition 5.1} There exists $Z_0\in\SZ$ with $D(Z_0)=0$.
\endproclaim

\sub{Proof} Our plan is to find a good initial point $Z^*\in\SZ$
and then follow the flow of $-\grad D$ from $Z^*$; our choice of
initial point will guarantee that the flow will lie along a curve
$\SY\subset\SZ$ along which $D(Z)\bigm|_\SY$ will have a special form.
Both the argument for the existence of a good initial point and the
argument that the negative gradient flow on the curve $\SY$ is only
critical at a point $Z$ with $D(Z)=0$ involve understanding how a
deformation of a zigzag affects extremal lengths on $\nez S$ and
$\sww S$, so we begin with that in subsection~5.2. In
subsection~5.3 we choose our good initial point $Z^*$, while in
subsection~5.4 we check that the negative gradient flow from $Z^*$
terminates at a reflexive symmetric zigzag. This will conclude
the proof of the main theorems.

\bpage
\sub{5.2. The tangent space to $\SZ$} In this subsection, we
compute a variational theory for zigzags appropriate for our
problem. In terms of our search for minimal surfaces, we recall
that the zigzags (and the resulting Euclidean geometry on
the domains $\nez\Om$ and $\sww\Om$) are
constructed to solve the period
problem for the \wei data: since we are left to show that
the domains $\nez\Om$ and $\sww\Om$
are conformally equivalent, we need a
formula for the variation of the extremal length (conformal
invariants) in terms of the periods.

More particularly, note again that what we are doing throughout
this paper is relating the Euclidean geometry of $\nez\Om$ (and $\sww\Om$,
respectively) with the conformal geometry of $\nez\Om$ (and $\sww\Om$, resp.)
The Euclidean geometry is designed to control the periods of the one-forms
$\nez\om$ (and $\sww\om$) and is restricted by the requirements that the
boundaries of $\nez\Om$ ($\sww\Om$, resp.) have alternating left and right
orthogonal turns, and that $\nez\Om$ and $\sww\Om$ are complementary domains of
a zigzag $Z$ in $\BC$. Of course, we are
interested in the conformal geometry of
these domains as that is the focus of the Main Theorem~B.

In terms of a variational theory, we are interested in deformations of a zigzag
through zigzags: thus, informally, the basic moves consist of shortening or
lengthening individual sides while maintaining the angles at the vertices.
These moves, of course, alter the conformal structure of the complementary
domains, and we need to calculate the effect on conformal invariants (in
particular, extremal length) of these alterations; those calculations 
involve the \tec theory described in section 2.3, and form the bulk of this
subsection. We list the approach below, in steps.

Step 1) We consider a self-diffeomorphism $f_\e$ of $\BC$ which takes a given
zigzag $Z_0$ to a new zigzag $Z_\e$: this is given explicitly in \S5.2.1,
formulas 5.1. (There will be two cases of this, which will in fact require two
different types of diffeomorphisms, which we label $f_\e$ and $f^*_\e$; they
are related via a symmetry, which will later benefit us through an important
cancellation.) These diffeomorphisms will be supported in a neighborhood of a
pair of edges; later in Step~4, we will consider the effect of contracting the
support onto increasingly smaller neighborhoods of those pair of edges.

Step 2) Infinitesimally, this deformation of zigzag results in infinitesimal
changes in the conformal structures of the complementary domains, and hence
tangent vectors to the \tec spaces of these domains. As described in the
opening of \S2.3, those tangent vectors are given by Beltrami differentials
$\nez{\dot\mu}$ (and $\sww{\dot\mu}$) on $\nez\Om$ (and $\sww\Om$) and it is
easy to compute $\nez{\dot\mu}$ (and $\sww{\dot\mu}$) in terms of $\(\f
d{d\e}f_\e\)_{\bar z}$ and $\(\f d{d\e}f^*_\e\)_{\bar z}$. This is done
explicitly in \S5.2.1, immediately after the explicit computations of $f_\e$
and $f^*_\e$.

Step 3) We apply those formulas for $\nez{\dot\mu}$ and $\sww{\dot\mu}$ to the
computation of the derivatives of extremal lengths (e.g. $\f
d{d\e}\bigm|_{\e=0}\nez E(k)$, in the notation of \S4). \tec theory (\S2.3)
provides that this can be accomplished through formula (2.1), which exhibits
gradient vectors $d\nez\ext(k)$ as meromorphic quadratic differentials
$\Phi^{\nes}_k$ (and $\Phi^{\sw}_k$) on the sphere $\nez S$ (and $\sww S$,
resp.) As we described in \S2.3, Gardiner [Gar]  gives a recipe for
constructing these differentials in terms of the homotopy classes of their
leaves. We describe these differentials in \S5.2.2.

Step 4) We have excellent control on these quadratic differentials along the
(lift of the) zigzag. Yet the formula (2.1) requires us to consider an integral
over the support of the Beltrami differentials $\nez{\dot\mu}$ and
$\sww{\dot\mu}$. It is convenient to take a limit of
$$
\int\phi^{\nes}(k)\nez{\dot\mu},
$$
and the corresponding $\sw$ integral, as the support of $\nez{\dot\mu}$ is
contracted towards a single pair of symmetric segments. We take this limit and
prove that it is both finite and non-zero in \S5.2.3: the limit
then clearly has a
sign which is immediately predictable based on which segments of the zigzag are
becoming longer or shorter, and how those segments meet the curve whose
extremal length we are measuring. The main difficulty in taking the limits of
these integrals is in controlling the appearance of some apparent
singularities: this difficulty vanishes once one invokes the symmetry condition
to observe that the apparent singularities cancel in pairs.

We begin our implementation of this outline with some notation.
Choose a zigzag $Z$; let
$I_k$ denote the segment of $Z$ connecting the points $P_k$ and
$P_{k+1}$. Our goal is to consider the effect on the conformal
geometries of $\nez S$ and $\sww S$ of a deformation of $\SZ$, where
$I_k$ (and $I_{-k-1}$, resp.) move into $\nez\Om$: one of the
adjacent sides $I_{k-1}$ and $I_{k+1}$ ($I_{-k-2}$ and $I_{-k}$,
resp.) is shortened and one is lengthened, and the rest of the
zigzag is unchanged.  (See Figure 2.)

\epsfxsize 4in
\centerline{\epsfbox{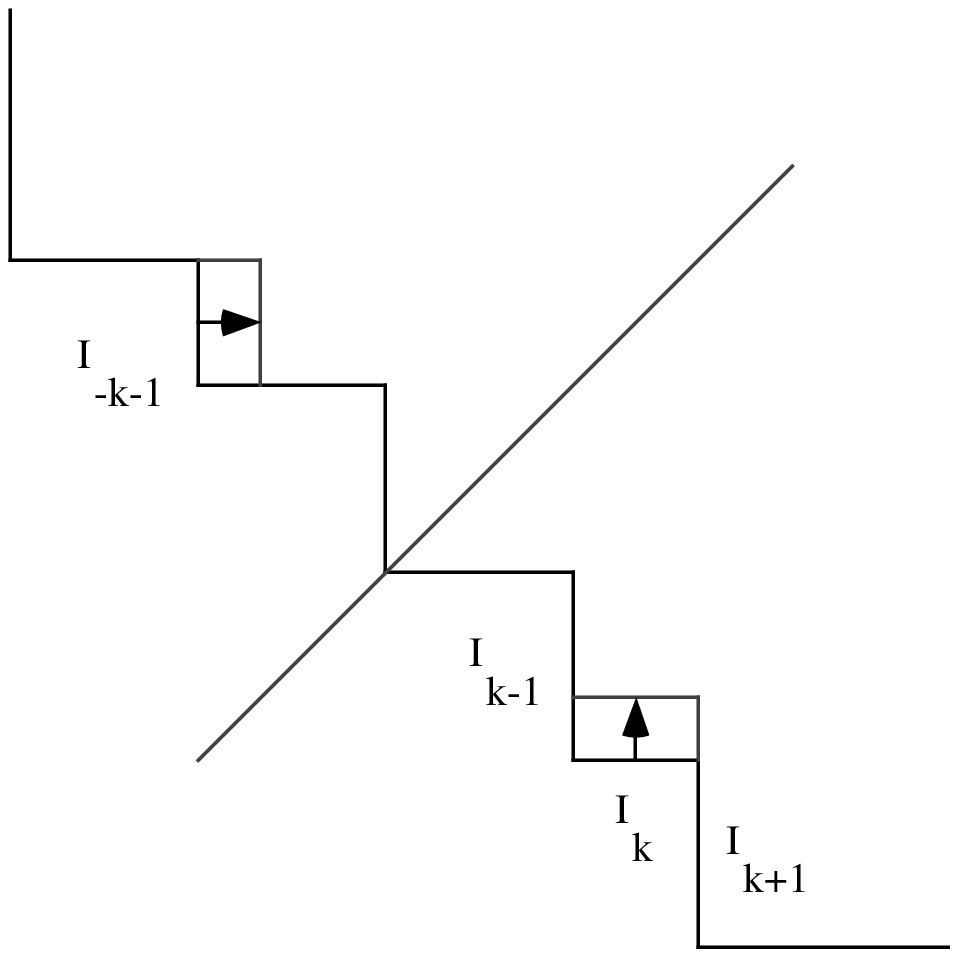}}
\centerline{Figure 2}

\sub{5.2.1} In this subsection, we treat steps 1 and 2 of the
above outline.
We begin by defining a family of maps $f_\e$ which move
$I_\e$ into $\nez\Om$; we presently treat the case that $I_k$ is
horizontal, deferring the vertical case until the next paragraph.

With no loss in generality, we may as well assume that $I_{k+1}$
is vertical; the more general case will just follow from obvious
changes in notations and signs. We consider local (conformal)
coordinates $z=x+iy$ centered on the midpoint of $I_k$ (i.e., the
horizontal segment abutting $I_k$ at the vertex $P_k$ of $I_k$
nearest to the line $\{y=x\}$.) In particular, suppose that $I_k$
is represented by the real interval $[-a,a]$, and define, for $b>0$
and $\d>0$ small, a local Lipschitz deformation $f_\e:\BC\to\BC$
$$
f_\e(x,y) = \cases
\(x,\e +\f{b-\e}b y\),   &\{-a\le x\le a, 0\le y\le b\}=R_1\\
\(x,\e +\f{b+\e}b y\),   &\{-a\le x\le a, -b\le y\le0\}=R_2\\
\(x,y+\f{\e+\f{b-\e}b y-y}\d(x+\d+a)\),   &\{-a-\d\le x\le-a,
0\le y\le b\} = R_3\\
\(x,y-\f{\e+\f{b-\e}b y-y}\d(x-\d-a)\),   &\{a\le x\le a+\d,
0\le y\le b\} = R_4\\
\(x,y+\f{\e+\f{b+\e}b y-y}\d(x+\d+a)\),   &\{-a-\d\le x\le-a,
-b\le y\le0\} = R_5\\
\(x,y-\f{\e+\f{b+\e}b y-y}\d(x-\d-a)\),   &\{a\le x\le a+\d,
-b\le y\le0\} = R_6\\
(x,y)   &\text{otherwise}
\endcases\tag"(5.1a)"
$$
where we have defined the regions $R_1,\dots,R_6$ within the
definition of $f_\e$. Also note that here $Z_0$ contains the arc
$\{(-a,y)\mid0\le y\le b\}\cup\{(x,0)\mid-a\le x\le a\}\cup
\{(a,y)\mid-b\le y\le0\}$.

\epsfxsize 4in
\centerline{\epsfbox{figure3.ps}}
\centerline{Figure 3}

Of course $f_\e$ differs from the identity only on a neighborhood
of $I_k$; so that $f_\e(Z_0)$ is a zigzag but no longer a symmetric
zigzag. We next modify $f_\e$ in a neighborhood of the reflected
(across the $y=x$ line) segment $I_{-k-1}$ in an analogous way with a
map $f^*_\e$ so that $f^*_\e\circ f_\e(Z)$ will be a symmetric
zigzag. (Here $f^*_\e$ is exactly a reflection of $f_\e$ if $k>1$.
In the case $k=1$, we require only a small adjustment for the fact
that $f_\e$ has changed one of the sides adjacent to the segment
$\ov{P_{-1}P_0}$, both segments of which lie in
$\supp(f^*_\e-\id)$.)

Our present conventions are that $I_k$ is horizontal; this forces
$I_{-k-1}$ to be vertical and we now write down $f^*_\e$ for such a
vertical segment; this is a straightforward extension of the
description of $f_\e$ for a horizontal side, but we present the
definition of $f^*_\e$ anyway, as we are crucially interested in
the signs of the terms. So set
$$
f^*_\e = \cases
\(\e +\f{b-\e}bx, y\),   &\{0\le x\le b, -a\le y\le a\}=R^*_1\\
\(\e +\f{b+\e}x, y\),   &\{-b\le x\le0, -a\le y\le a\}=R^*_2\\
\(x - \f{\e+\f{b-\e}b x-x}\d(y-\d-a), y\),   &\{0\le x\le b, a\le
y\le a+\d\}=R^*_3\\
\(x + \f{\e+\f{b-\e}b x-x}\d(y+\d+a), y\),   &\{0\le x\le b,
-a-\d\le y\le-a\}=R^*_4\\
\(x - \f{\e+\f{b+\e}b x-x}\d(y-\d-a), y\),   &\{-b\le x\le0, a\le
y\le a+\d\}=R^*_5\\
\(x + \f{\e+\f{b+\e}b -x}\d(y+\d+a), y\),   &\{-b\le x\le0, -a-\d\le
y\le-a\}=R^*_6\\
(x,y)   &\text{otherwise}
\endcases.\tag"(5.1b)"
$$
Note that under the reflection across the line $\{y=x\}$, the
regions $R_1$ and $R_2$ get taken to $R^*_1$ and $R^*_2$, but $R_4$
and $R_6$ get taken to $R^*_3$ and $R^*_5$, while $R_3$ and $R_5$
get taken to $R^*_4$ and $R^*_6$, respectively.

Let $\nu_\e = \f{\(f_\e\)_{\bar z}}{\(f_e\)_z}$
denote the Beltrami differential of $f_\e$, and set
$\dot\nu=\f d{d\e}\bigm|_{\e=0}\nu_\e$. Similarly, let $\nu^*_\e$
denote the Beltrami differential of $f^*_\e$, and set $\dot\nu^*=\f
d{d\e}\bigm|_{\e=0}\nu^*_\e$. Let $\dot\mu=\dot\nu+\dot\nu^*$. Now
$\dot\mu$ is a Beltrami differential supported in a bounded domain
in $\BC=\nez\Om\cup Z_0\cup\sww\Om$ around $Z_0$, so it restricts to a
pair $(\nez{\dot\mu},\sww{\dot\mu})$ of Beltrami differentials on the
pair of domains $(\nez\Om,\sww\Om)$. Thus, this pair of Beltrami
differentials lift to a pair $(\nez{\dot\mu},\sww{\dot\mu})$ on the
pair $(\nez S,\sww S)$ of punctured spheres, where we have
maintained the same notation for this lifted pair. But then, as a
pair of Beltrami differentials on $(\nez S,\sww S)\subset\tsym T$,
the pair $(\nez{\dot\mu},\sww{\dot\mu})$ represents a tangent
vector to $\SZ\subset\tsym T$ at $Z_0$. It is our plan to evaluate
$dD(\nez{\dot\mu},\sww{\dot\mu})$ to a precision sufficient to show
that $dD(\nez{\dot\mu},\sww{\dot\mu})<0$. To do this, we compute
$dExt([\g])$ for relevant classes of curves $[\g]$.

We begin by observing that it is easy to compute that
$\dot\nu=\(\f d{d\e}\bigm|_{\e=0}f_\e\)_{\bar z}$ evaluates near $I_k$ to
$$
\dot\nu = \cases
\f1{2b},   &z\in R_1\\
-\f1{2b},   &z\in R_2\\
\f1{2b}[x+\d+a]/\d+i\(1-y/b\)\f1{2d}=\f1{2b\d}(\bar z+\d+a+ib),
&z\in R_3\\
-\f1{2b}[x-\d-a]/\d-i\(1-y/b\)\f1{2\d}=\f1{2b\d}(-\bar z+\d+a-ib),
&z\in R_4\\
-\f1{2b}[x+\d+a]/\d+i\(1+y/b\)\f1{2\d}=\f1{2b\d}(-\bar z-\d-a+ib),
&z\in R_5\\
\f1{2b}[x-\d-a]/\d-i\(1+y/b\)\f1{2\d}=\f1{2b\d}(\bar z-\d-a-ib),
&z\in R_6\\
0  &z\notin\supp(f_\e-\id)
\endcases.\tag"(5.2a)"
$$
We further compute
$$
\dot\nu^* = \cases
-\f1{2b},   &R^*_1\\
\f1{2b},   &R^*_2\\
\f1{2b\d}(i\bar z-\d-a-bi)    &R^*_3\\
\f1{2b\d}(-i\bar z-\d-a+bi)   &R^*_4\\
\f1{2b\d}(-i\bar z+\d+a-bi)   &R^*_5\\
\f1{2b\d}(i\bar z+\d+a+bi)    &R^*_6
\endcases\tag"(5.2b)"
$$
Of course, this then defines the pair
$(\nez{\dot\mu},\sww{\dot\mu})$ by restriction to the appropriate
neighborhoods. In particular, $\nez{\dot\mu}$ is supported in the
(lifts of) the regions $R_1$, $R_4$, $R_6$, $R^*_1$, $R^*_3$ and
$R^*_5$ while $\sww{\dot\mu}$ is supported in the (lifts of) $R_2$,
$R_3$, $R_5$, $R^*_2$, $R^*_4$ and $R^*_6$.

\sub{5.2.2} We next consider the effect of the variation $\f
d{d\e}\bigm|_{\e=0}f_\e$ 
upon the conformal geometries of $\nez S$ and $\sww S$.
We compute the infinitesimal changes of some extremal lengths
induced by the variation $\f d{d\e}\bigm|_{\e=0}f_\e$.

For $[\g]$ a homotopy class of (a family of) simple closed curves,
the form $d\ext_{(\cd)}([\g])\in T^*_{(\cd)}\tsym T$ is given by an
element of $\qd^\sym(\cd)$. We describe some of these quadratic
differentials now; this is step 3 of the outline.

To begin, since the holomorphic quadratic differential
$\phi_k^{NE}=\phi_{\g_k}^{NE}=d\ext_{\nez S}([\g_k])$ is an element of
$\tsym\qd(\nez S)$, it is lifted from a holomorphic quadratic
differential $\psine k$ on ${\nez\Om}$ whose horizontal foliation has
nonsingular leaves either orthogonal to and connecting the segments
$I_{-k-2}$ and $I_{-k}$ or orthogonal to and connecting the
segments $I_k$ and $I_{k+2}$. (The foliation is parallel to the other
segments of $Z$, and the vertices where the foliation changes from
orthogonality to parallelism lift to points where the differential
$\phi_k^{NE}$ has a simple pole.)

Now the segments $I_k\in\ov{\nez\Om}$ corresponds under the map
$\phi:\nez\Om\to\sww\Om$ to the segment $I_{-k-2}\in\ov{\sww\Om}$;
similarly, $I_{-k-2}\in\ov{\nez\Om}$ corresponds to
$I_k\in\ov{\sww\Om}$. Thus $d\ext_{\sww S}\in\tsym\qd(\sww S)$
is lifted from a holomorphic quadratic differential $\psisw k$ whose
horizontal foliation has nonsingular leaves orthogonal to and
connecting the segments $I_{-k-2}$ and $I_{-k}$ and orthogonal to
and connecting $I_k$ and $I_{k+2}$, in an analogous way to  $\psine k$.
Now the foliations have characteristic local forms near the support of
the divisors of the differentials, and so the
foliations of $\phine$ (and $\phisw$) determine the divisors
of these differentials.
We collect this discussion, and its implications for the divisors, as

\bpage
\proclaim{Lemma 5.2} The horizontal foliations for $\psine k$ and
$\psisw k$ extend to a foliation of $\BC=\nez\Om\cup Z_0\cup\sww\Om$,
which is singular only at the vertices of $Z$. This foliation is
parallel to $Z$ except at $I_{-k-2}$, $I_{-k}$, $I_k$ and $I_{k+2}$,
where it is orthogonal. The differential $\phine$ (and $\phisw$)
have divisors
$$
\split
(\phine)  &= P^2_0P^2_\infty(P_{-k+1}P_{-k}P_{-k-1}P_{-k-2})^{-1}
(P_{k-1}P_kP_{k+1}P_{k+2})^{-1} = (\phisw)\ \text{ if }\ k>1\\
(\phi^{\nes}_1)  &= P^2_\infty(P_{-3}P_{-2}P_{-1}P_1P_2P_3)^{-1}
=(\phi^{\sw}_1)
\endsplit
$$
where $P_j$ refers to the lift of $P_j\in Z$ to $\nez S$ and
$\sww S$, respectively.
\endproclaim

\sub{5.2.3} Let $\phi^{\nes}$ denote a meromorphic quadratic
differential on $\nez S$ (symmetric about the lift of $\{y=x\}$)
lifted from a (holomorphic) quadratic differential $\psi^{\nes}$ on
(the open domain) $\nez\Om$; suppose that $\phi^{\nes}$ represents
the covector $d\ext_{\cd}([\g])$ in $T^*_{\nez S}\tsym T$ for some
class of curves $[\g]$. Formula (2.1) says that
$$
\f d{d\e}\biggm|_{\e=0}\ext_{S^\e_{\nes}}([\g]) = 2\ree\int_{\nez S}
\phi^{\nes}\dot\mu_{\nes} = 4\ree\int_{\nez\Om}
\psi^{\nes}\dot\nu_{\nes}\tag5.3
$$
where $S^\e_{\nes}$ is the punctured sphere obtained by
appropriately doubling $f_\e(\nez\Om)$.

The formula (5.3) is the basic variational formula that we will
use to estimate the changes in the conformal geometries of
$\nez\Om$ and $\sww\Om$ as we vary in $\SZ$. However,
in order to evaluate these integrals to a precision sufficient to
prove Proposition~5.1, we require a lemma. As background to the
lemma, note that $\dot\nu_{\nes}$ and $\sww{\dot\nu}$ depend upon
a choice of small constants $b>0$ and $\d>0$ describing the size of
the neighborhood of $I_k$ and $I_{-k-1}$ supporting $\dot\nu_{\nes}$
and $\sww{\dot\nu}$; on the other hand,  a hypothesis like the
foliation of $\psi^{\nes}$ is orthogonal to or parallel to $I_k$ and
$I_{-k-1}$ concerns the behavior of $\psi^{\nes}$ only at $I_k$ and
$I_{-k-1}$ (i.e., when $b=\d=0$). Thus, to use this information
about the foliations in evaluating the right hand sides of formula
(5.3), we need to have control on
$\ree\int_{\nez\Om}\psi^{\nes}\dot\nu_{\nes}$
as $b$ and $\d$ tend to zero. This is step 4 of the outline we gave
at the outset of section 5.2.

\bpage
\proclaim{Lemma 5.3} $\lim_{b\to0,
\d\to0}\ree\int_{\nez\Om}\psi^{\nes}\dot\nu_{\nes}$ exists
and is finite and non-vanishing.
Moreover, if $\psi^{\nes}$ has foliation either
orthogonal to or parallel to $I_k\cup I_{-k-1}$, then the sign of
the limit equals sgn~$(\psi^{\nes}\dot\nu_{\nes}(q))$ where $q$ is
a point on the interior of $I_k\cup I_{-k-1}$.
\endproclaim

\sub{Proof} On the interior of $I_k\cup I_{-k-1}$, the coefficients
of both $\psi^{\nes}$ and $\dot\nu_{\nes}$ have locally constant
sign; as we see from $\psi^{\nes}$ being either orthogonal or
parallel to $Z$ and symmetric, and from the form of
$\dot\nu_{\nes}$ in (5.2a) and (5.2b). We then easily check that
the sign of the product $\psi^{\nes}\dot\nu_{\nes}$ is constant on
the interior of $I_k\cup I_{-k-1}$, proving the final statement of
the lemma.

The only difficulty in seeing the existence of a finite limit as
$b+\d\to0$ is the possible presence of simple poles of
$\phi^{\nes}$ at the lifts of endpoints of $I_k\cup I_{-k-1}$.

To understand the singular behavior of $\psi^{\nes}$ near a vertex of
the zigzag, we begin by observing that on a preimage on
$\nez S$ of such a vertex, the quadratic differential has a simple
pole. Now let $\om$ be a local uniformizing parameter near the
preimage of the vertex on $\nez S$ and $\z$ a local uniformizing
parameter near the vertex of $Z$ on $\BC$. There are two cases to
consider, depending on whether the angle in $\nez\Om$ at the vertex
is (i) $3\pi/2$ or (ii) $\pi/2$. In the first case, the map from
$\nez\Om$ to a lift of $\nez\Om$ in $\nez S$ is given in coordinates
by $\om=(i\z)^{2/3}$, and in the second case by $\om=\z^2$. Thus, in
the first case we write $\phine=c\f{d\om^2}\om$ so that $\psine
k=-4/9c(i\z)^{-4/3}d\z^2$, and in the second case we write $\psine
k=4cd\z^2$; in both cases, the constant $c$ is real with sign
determined by the direction of the foliation.

With these expansions for $\psine k$ and $\psisw k$, we can compute
$d\ext([\g_k])[\dot\mu]$; of course, this quantity is given by
formula (2.1) as
$$
\align
d\ext_{\nez S}[(\g_k])[\dot\mu] &= 2\ree\int_{\nez S}\phi_k\dot\mu\\
&= 4\ree\int_{\nez\Om}\psi_k\dot\nu\\
&= 4\ree\(\int_{R_1\cup R'_1} + \int_{R_4\cup R^*_3} +
\int_{R_6\cup R^*_5}\) \psi_k\dot\nu.\tag5.4
\endalign
$$
Clearly, as $b+\d\to0$, as $|\dot\nu|=O\(\max\(\f1b,\f1\d\)\)$, we need
only concern ourselves with the contribution to the integrals of the
singularity at the vertices of $Z$ with angle $3\pi/2$.

To begin this analysis, recall that we have assumed that $I_k$ is
horizontal so that $Z$ has a vertex angle of $3\pi/2$ at $P_{k+1}$
and $P_{-k-1}$. It is convenient to rotate a neighborhood of
$I_{-k-1}$ through an angle of $-\pi/2$ so that the support of
$\dot\nu$ is a reflection of the support of $\dot\nu^*$ (see
equation (5.1)) through a vertical line. If the coordinates of
$\supp\dot\nu$ and $\supp\dot\nu^*$ are $z$ and $z^*$, respectively
(with $z(P_{k+1})=z^*(P_{-k-1})=0$), then the maps which lift
neighborhoods of $P_{k+1}$ and $P_{-k-1}$, respectively, to the
sphere $\nez S$ are given by
$$z\mapsto(iz)^{2/3}=\om \quad\hbox{ and}\quad z^*\mapsto(z^*)^{2/3}=\om^* .$$
Now the
poles on $\nez S$ have coefficients $c\f{d\om^2}\om$ and
$-c\f{d\om^{*2}}{\om^*}$,
respectively, so we find that when we pull back these poles from
$\nez S$ to $\nez\Om$, we have
$\psi^{\nes}(z)=-\f49c\ {dz^2}/{\om^2}$ while
$\psi^{\nes}(z^*)=-\f49c\ dz^2/(\om^*)^2$ in the coordinates $z$
and $z^*$ for $\supp\dot\nu$ and $\supp\dot\nu^*$, respectively.
But by tracing through the conformal maps $z\mapsto\om\mapsto\om^2$ on
$\supp\dot\nu$ and $z^*\mapsto\om^*\mapsto(\om^*)^2$, we see that
if $z^*$ is the reflection of $z$ through a line, then
$$\f1{(\om(z))^2}=1/{\ov{\om^*(z^*)^2}}$$
 so that the coefficients
$\psi^{\nes}(z)$ and $\psi^{\nes}(z^*)$ of
$\psi^{\nes}=\psi^{\nes}(z)dz^2$ near $P_{k+1}$ and of
$\psi^{\nes}(z^*)dz^{*2}$ near $P_{-k-1}$ satisfy
$\psi^{\nes}(z)=\ov{\psi^{\nes}(z^*)}$, at least for the singular
part of the coefficient.

On the other hand, we can also compute a relationship between the
Beltrami coefficients $\dot\nu(z)$ and $\dot\nu^*(z^*)$, in the
obvious notation, after we observe that $f_\e^*(z^*)=-\ov{f_\e(z)}$.
Differentiating, we find that
$$
\split
\dot\nu^*(z^*)  &= \dot f^*(z^*)_{\ov{z^*}}\\
&=-\ov{\dot f(z)}_{\ov{z^*}}\\
&=(\ov{\dot f(z)})_z\\
&= \ov{\dot f(z)_{\ov z}}\\
&= \ov{\dot\nu(z)}.
\endsplit
$$
Combining our computations of $\psi^{\nes}(z^*)$ and $\dot\nu(z^*)$
and using that the reflection $z\mapsto z^*$ reverses orientation,
we find that (in the coordinates $z^*=x^*+iy^*$ and $z=x+iy$) for
small neighborhoods $N_{\k}(P_{k+1})$ and $N_{\k}(P_{-k-1})$ of $P_{k+1}$
and $P_{-k-1}$ respectively,
$$
\split
\ree &\int\lm_{\supp\dot\nu\cap N_\k(P_{k+1})}\psi^{\nes}(z)
\dot\nu(z)dxdy +
\ree\int\lm_{\supp\dot\nu^*\cap N_\k(P_{-k-1})}\psi^{\nes}(z^*)
\dot\nu(z^*)dx^*dy^*\\
&= \ree\int\lm_{\supp\dot\nu\cap N_\k(P)}\psi^{\nes}(z)\dot\nu(z) -
\psi^{\nes}(z^*)\dot\nu(z^*)dxdy\\
&= \ree\int\lm_{\supp\dot\nu\cap N_\k}\psi^{\nes}(z)\dot\nu(z) -
\ov{[\psi^{\nes}(z)+O(1)]}\ \ov{\dot\nu(z)} dxdy\\
&= O(b+\d)
\endsplit
$$
the last part following from the singular coefficients summing to a
purely imaginary term while $\dot\nu=O\(\f1b+\f1\d\)$,
and the neighborhood has area $b\d$. This
concludes the proof of the lemma. \qed

\sub{5.3 A good initial point for the flow}
In this subsection, we seek a symmetric zigzag $Z^*$ of
genus~$p$ with the property that $E_{\nes}(k)=E_{\sw}(k)$ for
$k=2,\dots,p-1$. This will greatly simplify the height function
$D(Z)$ at $Z^*$. Our argument for the existence of $Z^*$ involves
the

\sub{Assumption 5.4} There is a reflexive symmetric zigzag of
genus~$p-1$.

\bpage
Since Enneper's surface can be represented by the zigzag of just
the positive $x$- and $y$- axes,
and we already have represented the Chen-Gackstatter surface of
genus one  by a zigzag in \S3, the initial step
of the inductive proof of this assumption is established.

The effect of the assumption is that on the codimension 1
face of $\p\ov\SZ$
consisting of zigzags with $P_{-1}=P_0=P_1$, there is a degenerate
zigzag $Z^*_0$ with $E_{\nes}(k)=E_{\sw}(k)$ for $k=2,\dots,p-1$.
Our goal in this subsection is the proof of

\bpage
\proclaim{Lemma 5.5} There is a family $Z^*_t\subset\SZ$ of
non-degenerate symmetric zigzags with limit point $Z^*_0$ where
each zigzag $Z^*_t$ satisfies $E_{\nes}(k)=E_{\sw}(k)$ for
$k=2,\dots,p-1$.
\endproclaim

\sub{Proof} We apply the implicit function theorem to a
neighborhood $V$ of $Z^*_0=(Z^*_0,0)$ in $\p\ov\SZ\x(-\e,\e)$,
where we will identify a neighborhood of $Z^*_0$ in $\ov\SZ$ with
a neighborhood $U$ of $Z^*_0$ in $\p\ov\SZ\x[0,\e)$. So our
argument will proceed in three steps: (i) we first define our
embedding of $U$ into $V$, (ii) we then show that the mapping
$\Phi:(Z,t)\mapsto(E_{\nes}(2)-E_{\sw}(2),\dots,
E_{\nes}(p-1)-E_{\sw}(p-1))$ is differentiable and then (iii)
finally we show that $d\Phi\bigm|_{\p\SZ\x\{0\}}$ is an isomorphism
onto $\BR^{p-2}$.
The first two steps are essentially formal, while the last step
involves most of the geometric background we have developed so far,
and is the key step in our approach to the gradient flow.

For our first step, normalize the zigzags in $U$
(as in section 4.1) so that $P_{-p}=i$
and $P_p=1$ and for $Z$ in $U$ near $Z^*_0$, and let
$(t(Z),a_2(Z),\dots,a_{p-1}(Z))$ denote the Euclidean lengths of
the segments $\<I_0,I_1,\dots,I_{p-2}\>$. Then for $Z\in U$, let $Z$
have coordinates $(\psi(Z),t)$ where $\psi(Z)\in\p\SZ$ has
normalized Euclidean lengths $a_2(Z),\dots,a_{p-1}(Z)$. It is easy
to see that $\psi:U\to\p\ov\SZ$ is a continuous and well-defined
map.

Next we verify that the map $\Phi$ is differentiable. We can
calculate the differential $D\Phi\bigm|_U$ by applying some of the
discussion of the previous subsection~5.2. For instance, the matrix
$D\Phi\bigm|_U(Z)$ can be calculated in terms of
$dE_{\nes}(k)\bigm|_Z[\dot\mu]$ where $\dot\mu$ corresponds to an
infinitesimal motion of an edge of $Z$, as in formula (5.3). Indeed
we see that as $\e\to0$, all of the derivatives
$d(E_{\nes}(k)-E_{\sw}(k))[\dot\mu]$ are bounded and converge: this
follows easily from observing that the quadratic differentials
$\phi_k$ are bounded and converge as $\e\to0$ and then applying
formulas (5.3) and (5.2). In fact, when $\supp\dot\mu$ meets
the lift of $I_0\cup I_{-1}$,
the same argument continues to hold, after we make one observation.
We observe that we can restrict our attention to where we are
sliding only the segments $I_1,\dots,I_{p-1}$ (and their
reflections) and not $I_0$ (and $I_{-1}$), as the tangent space is
$p-1$ dimensional; thus these derivatives are bounded as well.

This is all the differentiability we need for the relevant version
of the implicit function theorem.

Finally, we show that
$d\Phi\bigm|_{Z^*_0}:T_{Z^*_0}\SZ_{p-1}\to\BR^{p-2}$ is an
isomorphism. To see this it is sufficient to verify that this
linear map $d\Phi\bigm|_{Z^*_0}$ has no kernel. So let $v\in
T_{Z^*_0}\SZ_{p-1}$, so that
$$
v = \sum^{p-1}_{i=1} c_i\dot\nu_i
$$
where $c_i\in\BR$ and $\dot\nu_i$ refers to an infinitesimal
perturbation of $I_i$ and $I_{-i-1}$ into $\nez\Om$ (in the
notation for zigzags in $\SZ_p$: for $Z^*_0$, we have that $I_0$
and $I_{-1}$ have collapsed onto $P_0$.

Suppose, up to looking at $-v$ instead of $v$, that some $c_i>0$,
and let $\{i_j\}$ be the subset of the index set $\{1,\dots,p-1\}$
for which $c_{i_j}>0$. We consider the (non-empty) curve system
$\G$ of arcs connecting $I_{i_j}$ to the interval
$\ov{P_pP_\infty}$ and $I_{-i_j-1}$ to $\ov{P_\infty P_{-p}}$, let
$\vp_v$ denote a Jenkins-Strebel differential associated to this
curve system. By construction, sgn~$\vp_v$ is constant on the
interior of every interval, and $\vp_v>0$ on the interior of $I_i$
if and only if the index $i\in\{i_j\}$.

Thus, by Lemma 5.3 and formula (5.2), we see that both
$$
d\ext_\G(\nez\Om)[v] = \sum^{p-1}_{i=1}c_i\int_{I_i}\vp_v\dot\nu_i>0
\tag5.5a
$$
and from formulas (5.2) and (5.3), using that the horizontal (and vertical)
foliation(s) of $\vp_v$ extend to $\sww\Om$, that
$$
d\ext_\G(\sww\Om)[v] = \sum^{p-1}_{i=1}c_i\int_{I_i}\vp_v\dot\nu_i
<0.\tag5.5b
$$
Thus
$$
d\(\ext_\G(\nez\Om) - \ext_\G(\sww\Om)\)[v]>0.\tag5.6
$$
Now suppose that $d\Phi\bigm|_{Z^*_0}[v]=0$. Then by the
definition of $\Phi$, we would have that
$d\ext_\G(\nez\Om)[v]=d\ext(\sww\Om)[v]$, and since
$\<\ext_{\G_i}(\nez\Om)\>$ provides local coordinates in \tec space,
we would see that the \tec distance between $\nez\Om$ and $\sww\Om$
would infinitesimally vanish. But that would force
$\ext_{\nez\Om}(\G)-\ext_{\sww\Om}(\G)$ to vanish to first order, which
contradicts our computation (5.5).

We conclude that $d\Phi\bigm|_{Z^*_0}$ is an isomorphism, so that
the implicit function theorem yields the statement of the lemma.
\qed

\sub{5.4 The flow from the good initial point and the
proofs of the main results}
We now consider one of our ``good'' zigzags
$Z^*_t\in\SZ_p$ and use it as an initial point from which to flow
along $-\grad D(Z)$ to a reflexive symmetric zigzag.

Let $\SY\subset\SZ_p$ denote the set of genus~$p$ zigzags for which
$\ext_{\nez\Om}(\G_i)=\ext_{\sww\Om}(\G_i)$ for $i=2,\dots,p-1$. As
extremal length functions are in $C^1(\tsym T)$ by
Gardiner-Masur [GM],
we see that $\SY$ is
a piecewise $C^1$ submanifold of $\SZ_p$.
(We shall note momentarily that in our case, these
extremal length functions are real analytic.)  We consider the height
function $D$ restricted to the set $\SY$.

\bpage
\proclaim{Lemma 5.6} $D\bigm|_\SY$ is proper and is critical only at
points $Z\in\SY\subset\SZ$ for which $D(Z)=0$.
\endproclaim

\sub{Proof} The properness of $D\bigm|_\SY$ follows from the
properness of $D$, as shown in \S4. We next show that if $D(Z)\neq0$
for $Z\in\SY$, then there exists a tangent vector $v\in T_Z\SZ$ for
which $dD[v]<0$, and for which
$d\(\ext_{\nez\Om}(\G_i)-\ext_{\sww\Om}(\G_i)\)[v]=0$ so that $v$
lies tangent to a fragment of $\SY$ and infinitesimally reduces the
height $D$.

Indeed, we observe that
$$
D\bigm|_\SY = \[\(\exp\f1{E_{\nes(1)}}\) -
\(\exp\f1{E_{\sw(1)}}\)\]^2 + \[E_{\nes(1)} - E_{\sw(1)}\]^2
$$
as the other terms vanish.

Now, observe that $\SZ$ is a real analytic submanifold of the real
analytic product manifold $\tsym T\x\tsym T$, being defined in
terms of periods of a pair of holomorphic forms on the underlying
punctured spheres. Next we observe that $E_{\nes(j)}$ and
$E_{\sw(j)}$ are, for each $j$, real analytic functions on $\tsym
T$ with non-degenerate level sets. To see this note that the
extremal length functions correspond to just the energy of harmonic
maps from the punctured spheres to an interval, with the required
analyticity coming from Eells-Lemaire [EL], or directly from 4.5;
 the non-degeneracy
follows from Lemma~5.3,
if we apply any $\dot\nu$ of the form (5.2) to the zigzag
(this will be developed in more detail in the following paragraph).
Thus, the set $\SY$ acquires the structure of a real analytic
submanifold properly embedded in $\SZ$. As $\SY$ is one-dimensional
near $Z_0$, it is one-dimensional (with no boundary points)
everywhere.

Now, for $Z\in\SY\subset\SZ$ which is {\it not} a zero of $D$, we have
for any tangent vector $\dot\mu$ the formula
$$
\aligned
dD[\dot\mu] &= 2\(\exp\(\f1{\ext_{\nez S}}
([\g_1])\) - \exp\(\f1{\ext_{\sww S}}([\g_1])\)\)\\
&\bigg\{-\exp\(\f1{\ext_{\nez S}}([\g_1])\) \(\ext_{\nez
S}([\g_1])\)^{-2} d\ext_{\nez S}([\g_1])[\dot\mu]\\
&\qquad + \exp\(\f1{\ext_{\sww S}}([\g_1])\) \(\ext_{\sww
S}([\g_1])\)^{-2} d\ext_{\sww S}([\g_1])[\dot\mu]\bigg\}\\
&+\[\ext_{\nez S}([\g_1]) - \ext_{\sww S}([\g_1])\]
(d\ext_{\nez S}([\g_1])[\dot\mu] - d\ext_{\sww S}([\g_1])[\dot\mu])\\
\endaligned
$$
Then if we evaluate this expression for, say, $\dot\mu$ being given
by lifting an infinitesimal move of just one side as in formulas (5.1) and
(5.2), we
find by an argument similar to that for inequalities (5.5) and
(5.6) that $dD[\dot\mu]\neq0$. This concludes the proof of the
lemma. \qed

\bpage

\sub{Conclusion of the proof of Proposition 5.1} We argue
by induction. The union of the positive $x$-
and $y$-axes, is reflexive via the explicit map $z\mapsto iz^3$; this
verifies the statement of the proposition for genus $p=0$. There is
also a unique reflexive zigzag for genus $p=1$, after we make use of
the permissible normalization $P_1 =1$; we can verify that both
$\nez S$ and $\sww S$ are square tori, as in the first paragraph of
\S3.  In
general, once we are given a reflexive symmetric
zigzag of genus~$p-1$, we are
able to satisfy Assumption~5.4, so that Lemmas~5.5 and 5.6
guarantee the existence of a one-dimensional submanifold
$\SY\subset\SZ$ along which the height function is proper. A
minimum $Z_0$ for this height function is critical for $D$ along
$\SY$, and hence satisfies $D(Z_0)=0$ by Lemma~5.6. \qed

As the
discussion of subsection~4.3 shows that if $D(Z_0)=0$, then
$Z_0$ is reflexive, we conclude from Proposition 5.1 the

\proclaim{Main Theorem B} There exists a reflexive symmetric zigzag
of genus ~$p$ for $p\ge0$ which is isolated in $\SZ_p$.
\endproclaim

\sub{Proof of Main Theorem B} The local uniqueness
follows from inequality (5.6) and the argument
following it. \qed

Our main goal then follows.

\sub{Proof of Main Theorem A} By Main Theorem B, there exists a
symmetric reflexive zigzag of genus~$p$. By Theorem~3.3, and
Lemma~3.4 from this zigzag we can find \wei data for the required
minimal surface.

\sub{5.5 A remark on a different height function}
In this subsection, we try to give some context to
the methods we adopted in Sections~4 and~5 by considering an
alternate and perhaps more natural height function.

Define a new height function
$$
H(Z) = \{d_{\text{Teich}}(\nez\SR,\sww\SR\}^2.
$$
Certainly $H(Z)=0$ if and only if $Z$ is reflexive. Moreover, here
the gradient flow is much easier to work with, at least locally in
$\SZ_p$. Indeed, we observe

\proclaim{Lemma 5.7} $H(Z)=0$ if and only if $Z$ is critical for
$H(\cd)$ on $\SZ_p$.
\endproclaim

\sub{Sketch of Proof} Clearly $H(Z)\ge0$ and $H$ is $C^1$ (even real
analytic by the proof of Lemma~5.5) on
$\SZ_p\subset\tsym T\x\tsym T$, so if $H(Z)=0$, then $Z$ is critical
for $H$ on $\SZ_p$.

So suppose that $H(Z)\neq0$. Then $\nez\SR$ is not conformally
equivalent to $\sww\SR$, so we can look at the unique \tec map in
the homotopy class $[\phi:\nez S\to\sww S]$ (see the opening of
\S3).

From the construction of $\nez S$, $\sww S$ and $\phi$ from
$\nez\Om$, $\sww\Om$ and $\phi$, we can draw many conclusions about
the \tec differentials $q\in\qd^{\sym}(\nez S)$ and
$q'\in\qd^{\sym}(\sww S)$. For instance, the foliations are either
perpendicular or parallel to the image of $Z$ on $\nez S$ and $\sww
S$, and by the construction of the \tec maps, there are zeros in
the lift of an interval $I_k$ on $\nez S$ if and only if there is a
corresponding zero in the lift of $\phi(I_k)$ on $\sww S$.
Moreover, there is a simple pole at a puncture on $\nez S$ if and
only if there is a
simple pole at a corresponding puncture on $\sww S$.

Finally, we observe that $q$ has only simple poles and as the
foliation of $q$ is symmetric about both the images of $Z$ and the
line $\{y=x\}$, we see that there cannot be a simple pole at either
the lift of $P_0$ or $\infty$. Yet
by Riemann-Roch, as there are $4$ more poles than
zeros, counting multiplicity, there must be a pair of intervals
$I_k\cup I_{-k-1}$ whose lift has no zeros of $q$ (and whose image
under $\phi$ of lift has no zeros of $q'\in\qd^{\sym}(\sww S)$).

Observe next that Kerckhoff's formula (2.2) says that the
horizontal measured foliation $(\SF_q,\mu_q)$ of $q$ extremizes the
quotient on the right hand side of 2.2. However, consider a
deformation (5.1) on our zero-free intervals $I_k\cup I_{-k-1}$. By
Lemma~5.3 and formula (5.3) we see that either
$$
d\ext_{\nez S}(\mu_q)[\dot\nu] > 0\quad\text{and}\quad
d\ext_{\sww S}(\mu_q)[\dot\nu] < 0
$$
or
$$
d\ext_{\nez S}(\mu_q)[\dot\nu] < 0\quad\text{and}\quad
d\ext_{\sww S}(\mu_q)[\dot\nu] > 0.
$$
In either case, we see from Kerckhoff's formula that
$dH\bigm|_Z[\dot\mu]\neq0$, and so $Z$ is not critical. \qed

So why use our more complicated height function? The answer lies in
formula (4.2), which combined with Kerckhoff's formula (2.1) shows
that $H(\cd):\SZ_p\to\BR$ is not a proper function on $\SZ_p$.
Thus, the backwards gradient flow might flow to a reflexive
symmetric zigzag, or it may flow to $\p\ov\SZ_p$.
A better understanding of this height function $H$ on $\SZ$ would
be interesting both in its own right and if it would lead to a
new numerical algorithm for finding minimal surfaces experimentally.

\newpage

\sub{\S6. Extensions of the Method: The Karcher-Thayer Surfaces}

Thayer [T], following work of Karcher [K], defined \wei data
(depending upon unknown constants) for a family of surfaces
$M_{p,k}$ of genus~$p(k-1)$ with one Enneper-type end of winding
order $2k-1$. In this notation, the surfaces of genus~$p$ described
in Theorem~A are written $M_{p,2}$. Karcher [K] has solved the
period problem for the surfaces $M_{0,k}$ and $M_{1,k}$ for $k>2$
and Thayer has solved the period problem for $M_{2,k}$ for $k>2$. He
has also found numerical evidence supporting the solvability of the
period problem for $p\le34$, $k\le9$. Here we prove

\proclaim{Theorem C} There exists a minimally immersed surface
$M_{p,k}$ of genus~$p(k-1)$ with one Enneper-type end with winding
order $2k-1$.
\endproclaim

\sub{Proof} We argue in close analogy with our proof of
Theorem A ($k=2$). Let $\SZ_{p,k}$ denote the space of equivalence
classes of zigzags with $2p+1$ finite vertices and angles at the
vertices alternating between $\pi/k$ and $\f{2k-1}k\pi$. We double
the complementary domains $\nez\Om^k$ and $\sww\Om^k$ of $Z$ to
obtain cone-metric spheres with cone angles of alternating $2\pi/k$
and $2\pi/k(2k-1)$. We then take a $k$-fold cover of those spheres,
branched at the images of the vertices of the zigzag on the
cone-metric spheres, to obtain Riemann surfaces $\nez\SR^k$ and
$\sww\SR^k$. By pulling back the form $dz$ from $\nez\Om^k$ and
$\sww\Om^k$, we obtain, as in Section~3, a pair of forms $\a=gdh$
and $\b=g^{-1}dh$ on which we can base our \wei representation. As
before, we set $dh=d\pi$, where $\pi$ is the branched covering map
$\pi:\nez\SR^k\to\wh\BC$, so that $dh$ is exact; as before, the
periods of $\a$ and $\b$ are constructed to be conjugate, as soon
as the zigzag is reflexive, and the induced metric (2.0) is
regular at the lifts of the finite vertices of the zigzag.

To see that we can find a reflexive zigzag within $\SZ_{p,k}$, we
observe that
by the remark at the end of section 4.4,
the same real non-analyticity arguments of Section~4.4
carry over to the present case, once we replace the $\pi/2$ and
$3\pi/2$ angles with $\pi/k$ and $\pi/k(2k-1)$ angles.
All of the rest of the arguments of
Section~4 carry over without change and we conclude that the height
function $D(Z)$ is proper on $\SZ_{p,k}$.

For the gradient flow, we can write down deformations along the
zigzag analogous to formula (5.1)
(it is sufficient just to conjugate the maps in 5.1 by a map
which shears the original zigzag with vertex angles $\pi/2$ and
$3\pi/2$ to a zigzag with angles $\pi/k$ and $\f{2k-1}k\pi$)
and then check that the proof of
Lemma~5.3 continues to hold. The rest of the arguments in
Section~5 carry over unchanged to the present case. Thus we find a
reflexive symmetric zigzag in $\SZ_{p,k}$ via the proof of Theorem~B.
The present theorem then follows. \qed

\sub{Remark} Of course, the arguments in the last two paragraphs
of the proof of Theorem C apply equally well to zigzags of
arbitrary alternating angles $\th$ and $2\pi - \th$, so one might
well ask why we do not generalize the statement of Theorem C even
farther. The answer lies in that while the \tec theory of sections 4
and 5 extends to zigzags with non-orthogonal angles, the discussion
in section 3 of the transition from the zigzags to regular minimal
surfaces only extends to the zigzags described in the proof of
Theorem C.  For instance, we of course require a finitely branched cover
over a double of the zigzag in order to get a surface of finite genus,
so we must restrict our attention at the outset to zigzags with
rational angles.  However, if the smaller angle should be of the
form $\th=\f{r}s\pi$, we find that an $s$-sheeted cover of the double
of a zigzag would be forced to have an induced metric (2.0) which
was not regular at the lifts of the finite vertices.

\end
\newpage



\Refs
\refstyle{A}
\widestnumber\key{halloooo}

\ref
\key{ Blo}
\by D. Blo{\ss}
\book Elliptische Funktionen und vollst\"andige Minimalfl\"achen
\bookinfo PH.D. Thesis
\publ Freie Universit\"at
\publaddr Berlin
\yr 1989
\endref

\ref
\key CG    \by  C.C. Chen and F. Gackstatter
\paper Elliptische und Hyperelliptische Function und vollstandige
Minimal fl\"achen von Enneparschan Typ
\jour Math. Ann.  \vol259    \pages359--369 \yr1982
\endref

\ref
\key Cos    \by  C. Costa
\paper Example of a complete minimal immersion in $\BR^3$ of genus
one and three embedded ends
\jour Bull. Soc. Bras. Mat.  \vol15    \pages47--54 \yr1984
\endref

\ref
\key EL \by J. Eells and L. Lemaire
\paper Deformations of Metrics and Associated Harmonic Maps
\paperinfo Patodi Memorial Vol. Geometry and Analysis (Tata Inst.,
1981), 33--45
\endref



\ref
\key{Esp}
\by N. do Esp\'\i rito-Santo
\book Superf\'\i cies m\'\i nimas completas em $\Bbb R^3$ com fim de tipo
Enneper
\bookinfo PH.D. Thesis (1992)
\publaddr University of Niteroi, Brazil
(published as: \paper Complete Minimal Surfaces with type Enneper End
\jour Ann. Inst. Fourier (Grenoble) \vol44 \yr1994 \pages525--577)
\endref

\ref
\key{FLP}
\by A. Fathi, F. Laudenbach and V. Poenaru
\book Traveaux de Thurston sur les Surfaces
\bookinfo Asterisquu
\vol66-67
\yr1979
\publ Societ\'e Mathematique de France
\publaddr Paris
\endref



\ref
\key{Gack}
\by F. Gackstatter
\paper \"Uber die Dimension einer Minimalfl\"ache und zur Ungleichung von St.
Cohn-Vossen
\jour Arch. Rational Mech. Annal.
\vol 61(2)
\pages 141--152
\yr 1975
\endref

\ref
\key Gar \by F. Gardiner
\book Teichmuller Theory and Quadratic Differentials
\publ Wiley Interscience \publaddr New York \yr1987
\endref

\ref
\key GM \by F.P. Gardiner and H. Masur
\paper Extremal length Geometry of Teichm\"uller Space
\jour Complex Analysis and its Applications \yr1991 \vol16 \pages209--237
\endref

\ref
\key Ho-Me \by D. Hoffman and W.H. Meeks III
\paper Embedded Minimal Surfaces of Finite Topology
\jour Ann. of Math. \yr1990 \vol131 \pages1--34
\endref

\ref
\key{Ho-Ka}
\by D. Hoffman, H. Karcher
\paper Complete embedded minimal surfaces of finite total curvature
\jour Geometry V (R. Osserman, ed.) Encyclopaedia of Mathematical
Sciences 
\vol90
\publ Springer
\publaddr Berlin
\yr1997
\endref


\ref
\key  HM   \by J. Hubbard and H. Masur
\paper Quadratic Differentials and Foliations
\jour Acta Math.   \vol142   \yr1979   \pages221--274
\endref



\ref
\key J \by J.A. Jenkins
\paper On the Existence of certain general Extremal Metrics
\jour Ann. of Math. \yr1957 \vol66 \pages440--453
\endref

\ref
\key Ke \by S. Kerckhoff
\paper The asymptotic geometry of Teichm\"uller Space
\jour Topology \vol19 \yr1980 \pages23--41
\endref




\ref
\key{J-M}
\by L. Jorge, W.H. Meeks III
\paper The topology of complete minimal surfaces of finite total Gaussian
curvature
\jour Topology
\vol 22(2)
\pages 203--221       \yr 1983
\endref

\ref
\key{Kar}
\by H. Karcher
\paper Construction of minimal surfaces
\inbook Surveys in Geometry
\pages 1--96        \yr 1989
\publ University of Tokyo
\endref

\ref
\key{Laws}
\by H.B. Lawson, Jr.
\book Lectures on Minimal Submanifolds
\publ Publish or Perish Press
\publaddr Berkeley
\yr 1971
\endref

\ref
\key{Lop}
\by F.J. Lopez
\paper The classification of complete minimal surfaces with total
curvature greater than $-12\pi$
\jour Trans. Amer. Math. Soc.
\vol 334(1)   \pages 49--74
\yr 1992
\endref

\ref
\key Oht \by M. Ohtsuka
\book Dirichlet Problem, Extremal Length, and Prime Ends
\publ Van Nostrand Reinhold
\publaddr New York
\yr 1970
\endref

\ref
\key{ Oss1}
\by R. Osserman
\paper Global properties of minimal surfaces in $E^3$ and $E^n$
\jour Annals of Math.
\yr 1964
\vol 80(2)
\pages 340--364
\endref

\ref
\key{Oss2}
\by R. Osserman
\book A Survey on Minimal Surfaces
\bookinfo 2nd edition
\yr 1986
\publ Dover Publications
\publaddr New York
\endref

\ref
\key Rain \by E.D. Rainville
\book Intermediate Differential equations
\bookinfo 2nd edition
\publ The Macmillan Company
\publaddr New York
\endref

\ref
\key S \by K. Sato
\paper Existence proof of One-ended Minimal Surfaces with Finite
Total Curvature
\jour Tohoku Math. J. (2) \vol48 \yr1996 \pages229--246
\endref

\ref
\key Sch  \by  R. Schoen
\paper Uniqueness, Symmetry and Embeddedness of Minimal Surfaces
\jour J. Diff. Geometry  \vol18  \yr1983   \pages791--809
\endref

\ref
\key Str   \by  K. Strebel
\book Quadratic Differentials
\publ Springer   \publaddr Berlin  \yr1984
\endref

\ref
\key Tha  \by  E. Thayer
\paper Complete Minimal Surfaces in Euclidean $3$-Space
\paperinfo Univ. of Mass. Thesis, 1994
\endref

\ref
\key WW \by M. Weber and M. Wolf
\paper \tec Theory and Handle Addition for Minimal Surfaces
\paperinfo in preparation
\endref

\ref
\key Wo \by M. Wolf
\paper On Realizing Measured Foliations via Quadratic
Differentials of Harmonic Maps to $\BR$-Trees
\jour J. D'Analyse Math \vol68  \yr1996 \pages107--120
\endref


\endRefs
